\documentclass[a4paper,12pt]{amsart}
\usepackage{mathrsfs}
\usepackage{amsmath,amssymb,latexsym,amsfonts,amscd}

\title{Explicit birational geometry of threefolds of general type, II}
\author{Jungkai A. Chen and Meng Chen}
\address{\rm Department of Mathematics, National Taiwan University, Taipei,
106, Taiwan} \email{jkchen@math.ntu.edu.tw}

\address{\rm Institute of Mathematics, Fudan University,
Shanghai, 200433, People's Republic of China}
\email{mchen@fudan.edu.cn}

\thanks{The first author was partially supported by TIMS, NCTS/TPE
and National Science Council of Taiwan. The second author was
supported by National Outstanding Young Scientist Foundation
(\#10625103) and NNSFC Key project (\#10731030)}

\newcommand{\bQ}{{\mathbb Q}}
\newcommand{\bP}{{\mathbb P}}
\newcommand{\roundup}[1]{\lceil{#1}\rceil}
\newcommand{\rounddown}[1]{\lfloor{#1}\rfloor}

\newcommand\Vol{\text{\rm Vol}}

\newcommand\OO{{\mathcal{O}}}

\newtheorem{thm}{Theorem}[section]
\newtheorem{lem}[thm]{Lemma}
\newtheorem{cor}[thm]{Corollary}
\newtheorem{prop}[thm]{Proposition}

\newtheorem{op}[thm]{Open Problem}
\theoremstyle{definition}
\newtheorem{defn}[thm]{Definition}
\newtheorem{setup}[thm]{}

\newtheorem{exmp}[thm]{Example}
\newtheorem{conj}[thm]{Conjecture}
\newtheorem{rem}[thm]{Remark}
\theoremstyle{remark}

\newtheorem{THM}{\bf Theorem}

\begin{document}
\begin{abstract}
Let $V$ be a complex nonsingular projective 3-fold of general type.
We shall give a detailed classification up to baskets of
singularities on a minimal model of $V$. We show that the $m$-canonical map of $V$ is
birational for all $m\geq 73$ and that the canonical volume
$\Vol(V)\geq \frac{1}{2660}$. When $\chi(\OO_V)\leq 1$, our result
is $\Vol(V)\geq \frac{1}{420}$, which is optimal. Other effective
results are also included in the paper.
\end{abstract}
\maketitle
\pagestyle{myheadings} \markboth{\hfill J. A. Chen and M. Chen
\hfill}{\hfill Explicit birational geometry of threefolds\hfill}

\section{\bf Introduction}

This is the continuation of our previous paper \cite{Explicit_I},
which also serves as a guidance for the history of this topic. In
this note, we will extend our technique in \cite{Explicit_I}, while
improving other known methods, to systematically study the
birational geometry of 3-folds of general type.

Recall that we have already proved the following:

\begin{THM} (\cite[Theorem 1.1]{Explicit_I}) {\em Let $V$ be a nonsingular
projective 3-fold of general type. Then
\begin{itemize}
\item[(1)] $P_{12}>0$;

\item[(2)] $P_{m_0}\geq 2$ for some positive integer $m_0\leq 24$.
\end{itemize}}
\end{THM}

Our main theorems of this paper are as follows.

\begin{thm}\label{m} Let $V$ be a nonsingular
projective 3-fold of general type. Then
\begin{itemize}
\item[(1)] $P_{m}>0$ for all $m\geq 27$.

\item[(2)] $P_{24}\geq 2$ and $P_{m_0}\geq 2$ for some positive integer $m_0\leq 18$.

\item[(3)] $\varphi_m$ is birational for all $m\geq 73$.

\item[(4)] $\Vol(V)\geq \frac{1}{2660}$. Furthermore,
$\Vol(V)=\frac{1}{2660}$ if and only if $P_2=0$ and
 either
$\chi(\OO_V)=3$, $\mathscr{B}(X)=\{ 9\times \frac{1}{2}(1,-1,1),
\newline 2\times \frac{1}{7}(1,-1,3), \frac{1}{19}(1,-1,7), 3\times
\frac{1}{3}(1,-1,1), \frac{1}{10}(1,-1,3), \frac{1}{4}(1,-1,1),
\newline \frac{1}{5}(1,-1,1)\}$
 or $\chi(\OO_V)=2$,
$\mathscr{B}(X)=\{2\times \frac{1}{2}(1,-1,1), 2\times
\frac{1}{7}(1,-1,3),\newline 2\times \frac{1}{5}(1,-1,2),
\frac{1}{19}(1,-1,7), \frac{1}{4}(1,-1,1)\}$ where $\mathscr{B}(X)$
is the basket of singularities on a minimal model $X$ of $V$.
\end{itemize}
\end{thm}

\begin{thm}\label{1} Let $V$ be a nonsingular
projective 3-fold of general type with $\chi(\OO_V)\leq 1$. Then
\begin{itemize}
\item[(1)] $\varphi_m$ is birational for all $m\geq 40$.

\item[(2)] $\Vol(V)\geq \frac{1}{420}$. Furthermore,
$\Vol(V)=\frac{1}{420}$ if and only if the basket of
singularities (on any minimal model $X$ of $V$) is $\{3\times
\frac{1}{2}(1,-1,1), \frac{1}{7}(1,-1,3), \frac{1}{5}(1,-1,2),
\frac{1}{4}(1,-1,1),\frac{1}{6}(1,-1,1)\}$.
\end{itemize}
\end{thm}

Theorem \ref{1} (2) is optimal due to the following example:

\begin{exmp}\label{ex1} (\cite[p151, No.23]{C-R} ) The
canonical hypersurface $X_{46}\subset \bP(4,5,6,7,23)$ has 7
terminal quotient singularities and the canonical volume
$K^3_{X_{46}}=\frac{1}{420}$. One knows $\chi(\OO_{X_{46}})=1$ since
$p_g(X_{46})= q(X_{46})=h^2(\OO_{X_{46}})=0$. Furthermore, it is
known that $\varphi_{m}$ is birational for all $m\geq 27$, but
$\varphi_{26}$ is not birational.
\end{exmp}

\medskip

Throughout, we will frequently use those definitions, equalities
and inequalities about formal baskets in our previous paper (see
\cite[Sections 3,4]{Explicit_I}).

\section{\bf Technical preparation}

In this section, we set up some notions and principles evolved in
our detailed study. We shall prove some general results on
pluricanonical birationality and the lower bound of canonical
volume. Though the method has already appeared in several previous
works, the way of applying it is resultful to the effect that we
are able to treat various situations while proving our main
theorems.

\begin{setup}\label{reduc}{\bf Reduction to problems on minimal 3-folds}. Let $V$ be a
nonsingular projective 3-fold of general type. By the 3-dimensional
Minimal Model Program (see, for instance, \cite{KMM,K-M,Reid83}),
$V$ has a minimal model $X$ (with $K_X$ nef and admitting
$\bQ$-factorial terminal singularities). Denote by $K_X$ a canonical
divisor of $X$. A basic fact is that $\Vol(V)=K_X^3>0$.
{}From the view point of birational geometry, it suffices to prove
our main theorems only for minimal 3-folds $X$.
\end{setup}

\begin{defn} (1) The number $\rho_i=\rho_i(X)$ denotes the minimal positive integer
such that $P_{m}(X)>i$ for all $m\geq \rho_i$ where $i=0,1$.

(2) The number $\mu_i=\mu_i(X)$ denotes the minimal positive integer
with $P_{\mu_i}=P_{\mu_i}(X)>i$ where $i=0,1,2$.

(3) Denote by $\mathscr{B}(X)$ the basket of singularities on $X$
(according to Reid \cite{YPG}), and by $r(X)$ the Cartier index of
$X$.

By our definition, we see $\rho_0\leq \rho_1$ and $\mu_0\leq
\mu_1\leq \rho_1$. The existence of $\rho_1$ can be guaranteed by
Theorem 1.
\end{defn}

Now suppose we have $P_{m_0}\geq 2$ for certain positive integer
$m_0$. We may study the geometry of the rational map
$\varphi_{m_0}:=\Phi_{|m_0K_{X}|}$.


\begin{setup}\label{setup}{\bf Set up for $\varphi_{m_0}$.} We
study the $m_0$-canonical map of $X$:
$$\varphi_{m_0}:X\dashrightarrow \bP^{P_{m_0}-1}$$ which is
only a rational map. First of all we fix an effective Weil divisor
$K_{m_0}\sim m_0K_X$. By Hironaka's big theorem, we can take
successive blow-ups $\pi: X'\rightarrow X$ such that:
\begin{itemize}
\item [(i)] $X'$ is smooth; \item [(ii)] the movable part of
$|m_0K_{X'}|$ is base point free; \item [(iii)] the support of the
union of $\pi^*(K_{m_0})$ and the exceptional divisors is of simple
normal crossings.
\end{itemize}

Set $g_{m_0}:=\varphi_{m_0}\circ\pi$. Then $g_{m_0}$ is a morphism
by assumption. Let $X'\overset{f}\longrightarrow
\Gamma\overset{s}\longrightarrow W'$ be the Stein factorization of
$g_{m_0}$ with $W'$ the image of $X'$ through $g_{m_0}$. In
summary, we have the following commutative diagram:\medskip

\begin{picture}(50,80) \put(100,0){$X$} \put(100,60){$X'$}
\put(170,0){$W'$} \put(170,60){$\Gamma$}
\put(112,65){\vector(1,0){53}} \put(106,55){\vector(0,-1){41}}
\put(175,55){\vector(0,-1){43}} \put(114,58){\vector(1,-1){49}}
\multiput(112,2.6)(5,0){11}{-} \put(162,5){\vector(1,0){4}}
\put(133,70){$f$} \put(180,30){$s$} \put(95,30){$\pi$}
\put(130,-5){$\varphi_{m_0}$}\put(136,40){$g_{m_0}$}
\end{picture}
\bigskip

Recall that
$$\pi^*(K_X):=K_{X'}-\frac{1}{r(X)}E_{\pi}$$ with $E_{\pi}$ effective since $X$ is terminal.
So we always have
$$\roundup{m\pi^*(K_X)}\leq mK_{X'}$$ for any integer $m>0$. Denote by
$M_{m_0}$ the movable part of $|m_0K_{X'}|$. One has
$$m_0\pi^*(K_X)=M_{m_0}+E_{m_0}'$$
for an effective ${\mathbb Q}$-divisor $E_{m_0}'$. In total, since
$$h^0(X',\rounddown{m_0\pi^*(K_X)})=h^0(X',\roundup{m_0\pi^*(K_X)})=P_{m_0}(X')=P_{m_0}(X),$$
one has:
$$m_0K_{X'}=M_{m_0}+(E_{m_0}'+\frac{m_0}{r(X)}E_{\pi})$$ where
$E_{m_0}'+\frac{m_0}{r(X)}E_{\pi}$ is exactly the fixed part of
$|m_0K_{X'}|$.

If $\dim(\Gamma)\geq 2$, a general member $S$ of $|M_{m_0}|$ is a
nonsingular projective surface of general type by Bertini's theorem
and by the easy addition formula for Kodaira dimension.

If $\dim(\Gamma)=1$, a general fiber $S$ of $f$ is an irreducible
smooth projective surface of general type, still by the easy
addition formula for Kodaira dimension. We may write
$$M_{m_0}=\underset{i=1}{\overset{a_{m_0}}\sum}S_i\equiv
a_{m_0}S$$ where $S_i$ is a smooth fiber of $f$ for all $i$ and
$a_{m_0}\ge \text{min}\{2P_{m_0}-2, P_{m_0}+g(\Gamma)-1\}$, by
considering the degree of the divisor $f_*(M_{0})$ on $\Gamma$.
\end{setup}

\begin{defn} We call $S$ (in \ref{setup}) {\it a generic irreducible element
of} the linear system $|M_{m_0}|$. Denote by $\sigma:
S\longrightarrow S_0$ the blow-down onto the smooth minimal model
$S_0$. By abuse of concepts, we define {\it a generic irreducible
element} of an arbitrary movable linear system on any projective
variety in a similar way.
\end{defn}

\begin{defn}(1)
Define the positive integer $p=p(m_0)$ as follows:
$$p=\begin{cases}   1 &\text{if}\ \dim (\Gamma)\geq 2\\
 a_{m_0} &\text{if}\ \dim(\Gamma)=1.
 \end{cases}$$

(2) To simplify our statements, we say that the fibration $f$ is of
type III ( resp. II, I) if $ \dim \Gamma=3$ (resp. $2,1$).
According to our needs, we would like to classify type I into more
delicate ones:
$$f\  \text{is of type}\
\begin{cases} I_q & \text{if}\  g(\Gamma)>0,\\
            I_3 & \text{if}\  g(\Gamma)=0, \ P_{m_0} \ge 3,\\
            I_p& \text{if}\  g(\Gamma)=0,\ p_g(S) >0,\\
            I_n& \text{if}\  g(\Gamma)=0,\ p_g(S) =0.\\
\end{cases} $$
\end{defn}
\medskip

\begin{setup}\label{inv}{\bf Invariants of the fibration.} Let $V$
be a smooth projective 3-fold and $f:V\longrightarrow \Gamma$ a
fibration onto a nonsingular curve $\Gamma$. Leray spectral
sequence tells that:
$$E_2^{p,q}:=H^p(\Gamma,R^qf_*\omega_V)\Longrightarrow
E^n:=H^n(V,\omega_V).$$  By Serre duality and  \cite[Corollary 3.2,
Proposition 7.6]{Kol}, one has the torsion-freeness of the sheaves
$R^i f_* \omega_V$ and the following formulae:
$$h^2({\mathcal O}_V)=h^1(\Gamma,f_*\omega_V)+h^0(\Gamma,R^1f_*\omega_V),$$
$$q(V):=h^1({\mathcal O}_V)=g(\Gamma)+h^1(\Gamma,R^1f_*\omega_V).$$
\end{setup}

\begin{setup}\label{BP}{\bf Birationality principles.} Let $Y$ be a
nonsingular projective variety on which there are two divisors $D$
and $M$. Assume that $|M|$ is base point free. Take the Stein
factorization of $\Phi_{|M|}$: $Y\overset{f}\longrightarrow
W\longrightarrow \bP^{h^0(Y,M)-1}$ where $f$ is a fibration onto a
normal variety $W$. Then the rational map $\Phi_{|D+M|}$ is
birational onto its image if one of the following conditions is
satisfied:
\begin{itemize}
\item [(i)] (\cite[Lemma 2]{T}) $\dim\Phi_{|M|}(Y)\geq 2$,
$|D|\neq \emptyset$ and $\Phi_{|D+M|}|_S$ is birational for a
general member $S$ of $|M|$.

\item [(ii)] (\cite[\S2.1]{MPCPS})
$\dim\Phi_{|M|}(Y)=1$, $\Phi_{|D+M|}$ can separate different general
fibers of $f$ and $\Phi_{|D+M|}|_F$ is birational for a general
fiber $F$ of $f$.
\end{itemize}
\end{setup}

\begin{rem}\label{separate} For the condition \ref{BP} (ii), one knows
that $\Phi_{|D+M|}$ can separate different general fibers of $f$
whenever $\dim\Phi_{|M|}(Y)=1$, $W$ is a rational curve and $D$ is
an effective divisor. (In fact, since $|M|$ can separate different
fibers of $f$, so can $|D+M|$.)
\end{rem}

 Mostly, we will come across such a situation that a positive integer $m$, a base
point free linear system $|G|$ on $S$ and the linear system
$|M_{m_0}|$ simultaneously satisfy the following assumptions.

\begin{setup}\label{asum}{\bf Assumptions}. Denote by $C$ a generic irreducible element of $|G|$.
\begin{itemize}
\item [(1)] The linear system $|mK_{X'}|$
distinguishes different generic irreducible elements of $|M_{m_0}|$
(namely, $\Phi_{|mK_{X'}|}(S')\neq \Phi_{|mK_{X'}|}(S'')$ for two
different generic irreducible elements $S'$, $S''$ of $|M_{m_0}|$).

\item [(2)] The linear system $|mK_{X'}|_{|S}$ on $S$ (as a
sub-linear system of $|{mK_{X'}}_{|S}|$) distinguishes different
generic irreducible elements of $|G|$. {\small (Or sufficiently, the
complete linear system
$$|K_{S}
+\roundup{(m-1)\pi^*(K_X)-S-\frac{1}{p}E_{m_0}'}_{|S}|$$
distinguishes different generic irreducible elements of $|G|$.)}

\end{itemize}
\end{setup}

\begin{setup}\label{xi}{\bf A lower bound of $K^3$.}
We keep the same notation as above. Since $\pi^*K_X$ is nef and big,
there is a  rational number $\beta >0$ such that
$\pi^*(K_X)|_S-\beta C$ is numerically equivalent to an effective
$\bQ$-divisor on $S$.

We further define the following quantities:
\begin{eqnarray*}
&&\xi:=(\pi^*(K_X)\cdot C)_{X'};\\
&&\alpha:=(m-1-\frac{m_0}{p}-\frac{1}{\beta})\xi;\\
&&\alpha_0:=\roundup{\alpha}.
\end{eqnarray*}

\end{setup}

One has
$$K^3\geq \frac{p}{m_0}\pi^*(K_X)^2\cdot S\geq
\frac{p\beta}{m_0}(\pi^*(K_X)\cdot C)=\frac{p \beta}{m_0} \xi.
\eqno{(2.1)}$$

So it is essential to estimate the rational number
$\xi:=(\pi^*(K_X)\cdot C)_{X'}$ in order to obtain the lower bound
of $K^3$. We recall the following:

\begin{thm}\label{technical}(\cite[Theorem 3.2]{Chen-Zuo}) Keep the notation as above.  The inequality:
$$\xi\geq \frac{\deg(K_C)+\alpha_0}{m}$$
holds if one of the following conditions is satisfied:
\begin{itemize}
\item [(i)] $\alpha>1$; \item [(ii)] $\alpha>0$ and $C$ is an even
divisor, i.e. $C\sim 2H$ for a divisor $H$ on $S$.
\end{itemize}

Furthermore, under Assumptions \ref{asum} (1) and  (2), the map
$\varphi_m:=\Phi_{|mK_{X'}|}$ is birational onto its image if one
of the following conditions is satisfied:
\begin{itemize}
\item [(i)] $\alpha > 2$;
 \item [(ii)] $\alpha \geq 2$ and $C$ is
not a hyper-elliptic curve on $S$.
\end{itemize}
\end{thm}

\begin{rem} \label{weak} In particular
the inequality $\xi\geq \frac{\deg(K_C)+\alpha_0}{m}$ in Theorem
\ref{technical} implies
$$ \xi \ge \frac{\deg(K_C)}{1+
\frac{m_0}{p}+\frac{1}{\beta}}\eqno{(2.2)}$$ since, whenever $m$ is
big enough so that $\alpha>1$,
$$ m \xi \ge \deg(K_C)+ \alpha_0 \geq\deg(K_C)+
(m-1-\frac{m_0}{p}-\frac{1}{\beta}) \xi.$$
\end{rem}


As long as we have fixed a linear system $|G|$ on $S$, we are able
to prove the effective non-vanishing of plurigenera as follows.

\begin{prop} \label{non-v} Assume $P_{m_0}\geq 2$ for some positive integer $m_0$.
 Then $P_m(X)>1$ for all integers $m > 1+\frac{m_0}{p}+ \frac{1}{\beta}$.
In particular, $\rho_0\leq \rho_1\leq \rounddown{2+\frac{m_0}{p}+
\frac{1}{\beta}}$.
\end{prop}

\begin{proof}
Assume $m > 1+\frac{m_0}{p}+ \frac{1}{\beta}$. Keep the same
notation as in \ref{setup}. Put
$$\mathcal{ L}_m:=(m-1)\pi^*(K_X)-\frac{1}{p}E_{m_0}' . $$ Then we
have: $|K_{X'}+\roundup{\mathcal{L}_{m}}|\subset |mK_{X'}|.$ Noting
that $$\mathcal{L}_m-S\equiv (m-1-\frac{m_0}{p})\pi^*(K_X)|_S$$ is
nef and big, the  Kawamata-Viehweg vanishing theorem (\cite{Ka,V})
yields the surjective map
$$H^0(X', K_{X'}+\roundup{ \mathcal{L}_m}) \to H^0(S,
(K_{X'}+\roundup{ \mathcal{L}_m })_{|S}). \eqno(2.3)$$ Since $S$ is
a generic irreducible element of a free linear system, one has
$\roundup{*}|_S \geq \roundup{*_{|S} } $ for any divisor $*$ on
$X'$. It follows that
$$(K_{X'}+\roundup{ \mathcal{L}_m })_{|S} \geq {K_{X'}}_{|S}+\roundup{
{\mathcal{L}_m}_{|S} }\sim K_S+\roundup{(\mathcal{L}_m-S)_{|S}}.
\eqno(2.4)$$

Note that there is an effective ${\mathbb Q}$-divisor $\hat{H}$ on
$S$ such that $\frac{1}{\beta}\pi^*(K_X)|_{S}\equiv C+\hat{H}$. We
consider
$$\mathcal{D}_m:= (\mathcal{L}_m -S)_{|S}-\hat{H}$$ on $S$.
Then, by assumption, the divisor $\mathcal{D}_m-C\equiv (m-1-
\frac{m_0}{p}-\frac{1}{\beta})\pi^*(K_X)|_{S}$ is nef and big. Thus
the Kawamata-Viehweg vanishing theorem again gives the following
surjective map
$$H^0(S,
K_{S}+\roundup{\mathcal{D}_m}) \longrightarrow  H^0(C, K_C + D),
\eqno (2.5)$$ where $D:=\roundup{\mathcal{D}_m-C}|_C$ is a divisor
on $C$. Because $C$ is a generic irreducible element of a free
linear system, we have $D \ge \roundup{(\mathcal{D}_m-C)_{|C}}.$ A
simple calculation gives $$ \deg(D) \ge (\mathcal{D}_m-C)\cdot C =
(m-1- \frac{m_0}{p}-\frac{1}{\beta})\xi= \alpha>0.
 $$
Noting that $g(C)\geq 2$ since $S$ is of general type, Riemann-Roch
formula on $C$ gives $h^0(C, K_C+D)\geq 2$. {}Finally, surjective
maps (2.3), (2.5) and inequality (2.4) imply the statement.
\end{proof}


We need the following lemma while studying type $I_p$, $I_n$ and
$I_3$ cases.

\begin{lem}\label{>0} Let $S$ be a nonsingular projective surface
of general type. Denote by $\sigma:S\longrightarrow S_0$ the
blow-down onto its minimal model $S_0$. Let $Q$ be a $\bQ$-divisor
on $S$. Then $h^0(S,K_S+\roundup{Q})\geq 2$ under one of the
following conditions:
\begin{itemize}
\item[(i)] $p_g(S)>0$, $Q\equiv \sigma^*(K_{S_0})+Q_1$ for some
nef and big $\bQ$-divisor $Q_1$ on $S$;

\item[(ii)] $p_g(S)=0$, $Q\equiv 2\sigma^*(K_{S_0})+Q_2$ for some
nef and big $\bQ$-divisor $Q_2$ on $S$.
\end{itemize}
\end{lem}
\begin{proof} First of all $h^0(S, 2K_S)=h^0(S,2K_{S_0})>0$ by the
Riemann-Roch theorem on $S$, which is a surface of general type.
Fix an effective divisor $R_0\sim l\sigma^*(K_{S_0})$, where
$l=1,2$ in cases (i) and (ii) respectively. Then $R_0$ is nef and
big and $R_0$ is 1-connected by \cite[Lemma 2.6]{Maga}. The
Kawamata-Viehweg vanishing theorem says $H^1(S,
K_S+\roundup{Q}-R_0)=0$ which gives the surjective map:
$$H^0(S, K_S+\roundup{Q})\longrightarrow
H^0(R_0,K_{R_0}+G_{R_0})$$ where
$G_{R_0}:=(\roundup{Q}-R_0)_{|R_0}$ with $\deg(G_{R_0})\geq
(Q-R_0)R_0=Q_l\cdot R_0>0$. The 1-connectedness of $R_0$ allows us
to utilize the Riemann-Roch (see Chapter II, \cite{BPV}) as in the
usual way. Note that $S$ is of general type. So $K_{S_0}^2>0$ and
$\deg(K_{R_0})=2p_a(R_0)-2=(K_S+R_0)R_0\geq 2$. By the Riemann-Roch
theorem on the 1-connected curve $R_0$, we have
$$h^0(R_0,K_{R_0}+G_{R_0})\geq \deg(K_{R_0}+G_{R_0})+1-p_a(R_0)\geq p_a(R_0)
\geq 2.$$ Hence $h^0(S,K_S+\roundup{Q})\geq 2$.
\end{proof}

\begin{prop}\label{nonvanishing} Assume $P_{m_0}\geq 2$ for some positive integer $m_0$.
Then $P_m \ge 2$ for $m\geq h(m_0)$ under one of the following
situations:
\begin{itemize}
\item[(i)] $h(m_0)=2m_0+3$ when $f$ is of type $I_p$;

 \item[(ii)] $h(m_0)=3m_0+4$ when $f$ is of type $I_n$;

  \item[(iii)] $h(m_0)=\rounddown{\frac{3m_0}{2}}+4$ when $f$ is of type $I_3$.
\end{itemize}
In particular, $\rho_0 \le \rho_1 \le 2m_0+3, 3m_0+4,
\rounddown{\frac{3m_0}{2}}+4$, respectively.
\end{prop}

\begin{proof} Keep the same notation as in \ref{setup}.
When $f$ is of type I, we have $p=a_{m_0}$.  By \cite[Lemma
3.3]{Chen-Zuo}, there is a sequence of rational numbers
$\{\hat{\beta}_n\}$ with $\hat{\beta}_n\mapsto \frac{p}{m_0+p}\geq
\frac{1}{m_0+1}$ such that
$$\pi^*(K_X)_{|S}-\hat{\beta}_n\sigma^*(K_{S_0})\equiv H_n$$
for an effective $\bQ$-divisor $H_n$.

We consider
$$\mathcal{D}'_m:=(\mathcal{L}_m-S)_{|_S}-(m-1-\frac{m_0}{p}) H_n
\equiv (m-1-\frac{m_0}{p})\hat{\beta}_n \sigma^*(K_{S_0}).$$

If, for $m>0$, $h^0(S, K_S+\roundup{\mathcal{D}'_m})\geq 2$, then
$h^0(S, K_S+\roundup{(\mathcal{L}_m-S)_{|_S}})\geq 2$. It follows
then $P_m \geq 2$ by surjective map (2.3) and inequality (2.4). We
can choose $h(m_0)$ according to the type of $f$.

When $f$ is of type $I_p$, we can pick a big number $n$ so that
$\hat{\beta}_n\geq \frac{1}{m_0+1}-\delta$ for some $0 < \delta \ll
1$. For $m \ge 2m_0+3$, we see $(m-1-\frac{m_0}{p})\hat{\beta}_n
>1$. By Lemma \ref{>0} and since $p_g(S)>0$, we know $h^0(S,
K_S+\roundup{\mathcal{D}'_m})\geq 2$. Thus we may take
$h(m_0)=2m_0+3$.

When $f$ is of type $I_n$, we still take a big number $n$ so that
$\hat{\beta}_n\geq \frac{1}{m_0+1}-\delta$ for some $0 < \delta \ll
1$. But, for $m \ge 3m_0+4$, we have
$(m-1-\frac{m_0}{p})\hat{\beta}_n
>2$. By Lemma \ref{>0} again, we see $h^0(S,
K_S+\roundup{\mathcal{D}'_m})\geq 2$. Thus we may take $h(m_0)=
3m_0+4$.

Finally when $f$ is of type $I_3$, we have $p\geq 2$. One may take a
big number $n$ so that $\hat{\beta}_n \geq \frac{2}{m_0+2}-\delta$
for some $0 < \delta \ll 1$. For $m\geq
\rounddown{\frac{3m_0}{2}}+4$, we have
$(m-1-\frac{m_0}{p})\hat{\beta}_n
>2$. Lemma \ref{>0} implies $h^0(S,
K_S+\roundup{\mathcal{D}'_m})\geq 2$. Thus we may take
$h(m_0)=\rounddown{\frac{3m_0}{2}}+4$. This completes the proof.
\end{proof}

\begin{lem}\label{a(1)} Assume $P_{m_0}(X)\geq 2$ for some positive integer $m_0$. Keep the same
notation as in \ref{setup}. Then, for $m\geq \rho_0+m_0$,
Assumptions \ref{asum} (1) is satisfied if $f$ is of type $III, II,
I_3, I_p$ or $I_n$.
\end{lem}

\begin{proof} Let $t>0$ be an integer. We consider the linear system
$|K_{X'}+\roundup{t\pi^*(K_X)}+M_{m_0}|\subset |(m_0+t+1)K_{X'}|$.
Since $K_{X'}+\roundup{t\pi^*(K_X)}\geq (t+1)\pi^*(K_X)$, we see
that $K_{X'}+\roundup{t\pi^*(K_X)}$ is effective whenever $t+1\geq
\rho_0$.

When $f$ is of type $I_3, I_p$ or $I_n$, we necessarily have
$g(\Gamma)=0$. Thus, by \cite[Lemma 2]{T} and Remark
\ref{separate}, the linear system
$|K_{X'}+\roundup{t\pi^*(K_X)}+M_{m_0}|$ can separate different
generic irreducible elements $S$ of $|M_{m_0}|$.
\end{proof}

\begin{lem}\label{S2} Let $T$ be a nonsingular projective surface
of general type on which there is a base point free linear system
$|G|$. Let $Q$ be an arbitrary $\bQ$-divisor on $T$. Then the
linear system $|K_T+\roundup{Q}+G|$ can distinguish different
generic irreducible elements of $|G|$ under one of the following
conditions:
\begin{itemize}
\item[(i)] $K_T+\roundup{Q}$ is effective and $|G|$ is not
composed with an irrational pencil of curves;

\item[(ii)] $Q$ is nef and big and $|G|$ is composed with an
irreducible pencil of curves.
\end{itemize}\end{lem}
\begin{proof} Statement (i) follows from \cite[Lemma 2]{T} and
Remark \ref{separate}.

For statement (ii), we pick up a generic irreducible element $C$
of $|G|$. Then $G\equiv sC$ where $s\geq 2$ and $C^2=0$. Let $C'$
be another generic irreducible element. The Kawamata-Viehweg
vanishing theorem gives the surjective map:
$$H^0(T,K_T+\roundup{Q}+G)\longrightarrow H^0(C,K_C+D)\oplus
H^0(C',K_{C'}+D')$$ where $D:=(\roundup{Q}+G-C)|_C$ and
$D':=(\roundup{Q}+G-C')|_{C'}$ with $\deg(D)>0$, $\deg(D')>0$.
Since $T$ is of general type, both $C$ and $C'$ are curves of
genus $\geq 2$. Thus $h^0(C,K_C+D)=h^0(C',K_{C'}+D')>1$. Thus
$|K_T+\roundup{Q}+G|$ can distinguish $C$ and $C'$.
\end{proof}

\begin{lem}\label{a(2)} Assume $P_{m_0}(X)\geq 2$ for some positive integer
$m_0$. Keep the same notation as in \ref{setup}. Take $G:=S|_S$ for
a generic irreducible element $S$ of $|M_{m_0}|$. Then Assumptions
\ref{asum} (2) is satisfied under one of the following situations:
\begin{itemize}
\item[(i)] $f$ is of type $III$ and $m\geq \rho_0+m_0$.

\item[(ii)] $f$ is of type $II$ and $m\geq {\max}\{\rho_0+m_0,
2m_0+2\}$.
\end{itemize}
\end{lem}
\begin{proof}
Since $$\begin{array}{rrl}
&& K_S+\roundup{(m-1)\pi^*(K_X)-S-\frac{1}{p}E_{m_0}'}_{|S}\\
&\geq & K_S+(m-1)\pi^*(K_X)_{|S}-(S+E_{m_0}')_{|S}\\
&=& K_S+(m-m_0-1)\pi^*(K_X)_{|S}\\
&\geq& (m-m_0)\pi^*(K_X)_{|S}+G
\end{array}$$
and
$$\begin{array}{lll}
&&K_S+(m-m_0-1)\pi^*(K_X)_{|S}\\
&\geq& K_S+(m-2m_0-1)\pi^*(K_X)_{|S}+G,
\end{array}$$
Lemma \ref{S2} implies that
$|K_S+\roundup{(m-1)\pi^*(K_X)-S-\frac{1}{p}E_{m_0}'}_{|S}|$ can
distinguish different generic irreducible elements of $|G|$
respectively. Note that, if $f$ is of type $III$, $|G|$ is not
composed with a pencil of curves. We are done.
\end{proof}

Under the condition $P_{m_0}\geq 2$, we study the pluricanonical map
$\varphi_m$ according to the type of $f$.

\begin{setup} {\bf Type $III$.}\label{III}

When $f$ is of type $III$, we have $p=1$ by definition. In this
case, $S\sim M_{m_0}$ and $|S|$ gives a generically finite
morphism. We take $G:=S|_S$. Then $|G|$ is base point free and
$\varphi_{|G|}$ gives a generically finite map. So a generic
irreducible element $C\sim G$ is a smooth curve.

If $\varphi_{|G|}$ gives a birational map, then $\dim
\varphi_{|G|}(C)=1$ for a general member $C$. The Riemann-Roch and
Clifford's theorem on $C$ says $C^2=G\cdot C\geq 2$. If
$\varphi_{|G|}$ gives a generically finite map of degree $\geq 2$,
since $h^0(S,G)\geq h^0(X',S)-1\geq 3$, one gets $C^2\geq
2(h^0(S,G)-2)\geq 2$. Anyway we have $C^2\geq 2$. So
$\deg(K_C)=(K_S+C)\cdot C>2C^2\geq 4$. We see $\deg(K_C)\geq 6$
since it is a even number.

One may take $\beta=\frac{1}{m_0}$ since $m_0\pi^*(K_X)|_S\geq C$.

Now inequality (2.2) gives $\xi\geq \frac{6}{2m_0+1}$. Take
$m=3m_0+2$. Then $\alpha=(m-2m_0-1)\xi>3$. So, by Theorem
\ref{technical}, $\xi\geq \frac{10}{3m_0+2}$. It follows from
inequality (2.1) that $K^3\geq \frac{10}{(3m_0+2)m_0^2}$.

We now consider the non-vanishing of plurigenera.
 By Proposition \ref{non-v}, we have $P_m \geq 2$ for all $ m
> 2m_0+1$. Now, if $m=2m_0+1$, the surjective map (2.3) and inequality (2.4) lead us
to compute $h^0(S, K_S+ \roundup{m_0 \pi^*K_X|_S})$. Let $L$ be a
generic irreducible element in $|S|_S|$. Then $L$ is effective and
nef. Since $h^2(K_S+L)=0$, one has $h^0(S,
K_S+L)\geq\chi(S,K_S+L)=\frac{1}{2}(K_S\cdot
L+L^2)+\chi(\OO_S)\geq 2$ by Riemann-Roch theorem. Hence
$P_{2m_0+1}\geq 2$. Also, $P_{2m_0} \ge P_{m_0} \geq 2$.
Therefore, we have $P_m >1$ for all $m \ge 2m_0$. In particular,
$\rho_0 \le \rho_1 \le 2m_0$.

By Lemmas \ref{a(1)}, \ref{a(2)}, Assumptions \ref{asum} (1), (2)
are satisfied if  $m \ge 3m_0$. Now $\alpha=(m-2m_0-1) \xi \ge
(m-2m_0-1)\frac{10}{3m_0+2}$. One sees that $\alpha>2$ if $m >
\frac{13m_0+7}{5}$. Hence $\varphi_m$ is birational if $$m >\max\{
3m_0-1, \frac{13m_0+7}{5}\}.$$

We conclude the following:
\begin{thm}\label{tIII} Assume $P_{m_0}(X) \ge 2$ for some positive integer
$m_0$. If the induced map $f$ is of type $III$. Then
\begin{enumerate}
\item $\rho_0 \le \rho_1 \le 2m_0$.

\item $K^3\geq \frac{10}{(3m_0+2)m_0^2}$.

\item $\varphi_m$ is birational if $m>\max\{ 3m_0-1,
\frac{13m_0+7}{5}\}$.
\end{enumerate}
\end{thm}
\end{setup}

\begin{setup}{\bf Type $II$.}\label{II}

 When $f$ is of type $II$, we see that $S \sim M_{m_0}$. Take
$|G|:=|S|_S|$, which is, clearly, composed with a pencil of curves.

Since a generic irreducible element $C$ of $|G|$ is a smooth curve
of genus $\geq 2$, we have $\deg(K_C)\geq 2$. Furthermore we have
$h^0(S,G)\geq h^0(X',S)-1\geq 2$. So $G\equiv \widetilde{a} C$ where
$\widetilde{a}\geq h^0(S,G)-1\geq 1$. This means that
$m_0\pi^*(K_X)|_S\geq S|_S\geq_{\text{num}} C$. So we may take
$\beta = \frac{1}{m_0}$.

Now inequality (2.2) gives $\xi\geq \frac{2}{2m_0+1}$. Take
$m=3m_0+2$. Then $\alpha>1$. One gets $\xi\geq \frac{4}{3m_0+2}$ by
Theorem \ref{technical}. So inequality (2.1) implies: $K^3\geq
\frac{4}{(3m_0+2)m_0^2}$.

Exactly the same proof as in Type $III$ shows that $ \rho_0 \le
\rho_1 \le 2m_0$.

By Lemmas \ref{a(1)}, \ref{a(2)}, Assumptions \ref{asum} (1),(2)
are satisfied if  $m \ge 3m_0$. Now $\alpha=(m-2m_0-1) \xi \ge
(m-2m_0-1)\frac{4}{3m_0+2}$. One sees that $\alpha>2$ if $m >
\frac{7m_0+4}{2}$. Since $\frac{7m_0+4}{2}>3m_0$, $\varphi_m$ is
birational if $m >\frac{7m_0+4}{2}$.

We conclude the following:
\begin{thm}\label{tII} Assume $P_{m_0}(X) \ge 2$ for some positive integer
$m_0$. If the induced map $f$ is of type $II$. Then
\begin{enumerate}
\item $\rho_0 \le \rho_1 \le 2m_0$.

\item $K^3\geq \frac{4}{(3m_0+2)m_0^2}$.

\item $\varphi_m$ is birational if $m >
\frac{7m_0+4}{2}.$
\end{enumerate}
\end{thm}
\end{setup}

\begin{setup} \label{Iq}{\bf Type $I_q$.}

Since $g(\Gamma)>0$, one sees $q(X) >0$ and hence $X$ is irregular.
This case is particularly well-behaved. It's known that $\varphi_m$
is birational for all $m \ge 7$ (see \cite{JC-H}). Also $K_X^3 \ge
\frac{1}{22}$ (see \cite{JLMS}).
\end{setup}

\begin{setup} \label{Ip}{\bf Type $I_p$.}

We have an induced fibration $f:X'\longrightarrow \Gamma$ with
$g(\Gamma)=0$. By definition, $p=a_{m_0}\geq 1$. By assumption,
$p_g(S)>0$ for a general fiber $S$ of $f$. We take
$G:=2\sigma^*(K_{S_0})$. Then one knows that $|G|$ is base point
free (see \cite[Theorem 3.1]{Cili}). Thus $|G|$ is not composed with
a pencil and a generic irreducible element $C$ is smooth. By
\cite[Lemma 3.3]{Chen-Zuo}, we can find
 a sequence of rational numbers $\{{\beta}_n\}$ with ${\beta}_n\mapsto
\frac{p}{m_0+p}$ such that $\pi^*(K_X)_{|S}-
\frac{{\beta}_n}{2}C\equiv H_n$ for effective
$\bQ$-divisors $H_n$. We may assume that  $\beta \ge \frac{1}{2(m_0+1)}-\delta$ for some $0<\delta \ll 1$.

Since
$C\sim 2\sigma^*(K_{S_0})$,
$$\deg(K_C)=(K_S+C)\cdot C\geq (\pi^*(K_X)|_S+C)\cdot C>C^2\geq 4.$$
Since $\deg(K_C)$ is even, we see $\deg(K_C)\geq 6$.

Now inequality (2.2) gives $\xi\geq \frac{6}{3m_0+3}$. Take
$m=4m_0+5$. Then $\alpha=(m-1-m_0-\frac{1}{\beta})\xi >2$ and
Theorem \ref{technical} gives $\xi\geq \frac{9}{4m_0+5}$. So, by
inequality (2.1), one gets $K^3\geq \frac{9}{2m_0(m_0+1)(4m_0+5)}$.

 Note that
$$\begin{array}{rrl}
&& K_S+\roundup{(m-1)\pi^*(K_X)-S-\frac{1}{p}E_{m_0}'}_{|S}\\
&\geq & K_S+\roundup{(m-m_0-1)\pi^*(K_X)_{|S}}\\
&\geq& K_S+\roundup{(m-m_0-1)\pi^*(K_X)_{|S}-\frac{3}{\beta_n}H_n}\\
&=&K_S+\roundup{Q_1}+\sigma^*(K_{S_0})+C
\end{array}  \eqno{(2.6)}$$
where
$Q_1:=(m-m_0-1)\pi^*(K_X)_{|S}-C-\sigma^*(K_{S_0})-\frac{3}{\beta_n}H_n\equiv
(m-m_0-1-\frac{3}{\beta_n})\pi^*(K_X)|_S$ is nef and big whenever
$m\geq 4m_0+5$. By Lemma \ref{>0} (i),
$K_S+\roundup{Q_1}+\sigma^*(K_{S_0})$ is effective. Thus, according
to \cite[Lemma 2]{T}, Assumptions \ref{asum} (2) is satisfied for
$m\geq 4m_0+5$. Since Proposition \ref{nonvanishing} (ii) implies
$\rho_0\leq 2m_0+3$, Lemma \ref{a(1)} (ii) tells that Assumptions
\ref{asum} (1) is satisfied as long as $m\geq 3m_0+3$. Take $m\geq
4m_0+5$. Then $\alpha\geq (m-3m_0-3)\xi\geq
\frac{2m_0+4}{m_0+1}>2$. So Theorem \ref{technical} implies that
$\varphi_m$ is birational for all $m\geq 4m_0+5$.

We thus summarize:
\begin{thm}\label{tIp} Assume $P_{m_0}(X) \ge 2$ for some positive integer
$m_0$. If the induced map $f$ is of type $I_p$. Then
\begin{enumerate}
\item $\rho_0 \le \rho_1 \le 2m_0+3$.

\item $K^3\geq \frac{9}{2m_0(m_0+1)(4m_0+5)}$.

\item $\varphi_m$ is birational if $m \ge 4m_0+5$.
\end{enumerate}
\end{thm}
\end{setup}

\begin{setup} \label{In}{\bf Type $I_n$.}

Similar to type $I_p$ case, we have  $p\geq 1$.  We take
$|G|:=|4\sigma^*(K_{S_0})|$ which is base point free by a
well-known result in \cite{Bom}.  Thus $|G|$ is not composed with a
pencil and a generic irreducible element $C$ is smooth. Similarly,
we can find
 a sequence of rational numbers $\{{\beta}_n\}$ with ${\beta}_n\mapsto
\frac{p}{m_0+p}$ such that $\pi^*(K_X)_{|S}-
\frac{{\beta}_n}{4}C\equiv H_n$ for effective $\bQ$-divisors $H_n$.
We may assume that  $\beta \ge \frac{1}{4(m_0+1)}-\delta$ for some
$0<\delta \ll 1$.

Since $\deg(K_C)>16\sigma^*(K_{S_0})^2\geq 16$ and $\deg(K_C)$ is
even, inequality (2.2) gives $\xi\geq \frac{18}{5m_0+5}$. Take
$m=6m_0+6$. Then
$\alpha=(m-1-m_0-\frac{1}{\beta})\xi=\frac{18}{5}>3$ and Theorem
\ref{technical} gives $\xi\geq \frac{11}{3m_0+3}$. So, by
inequality (2.1), one gets $K^3\geq \frac{11}{12m_0(m_0+1)^2}$.

By Proposition \ref{nonvanishing}, we have $P_{m} \ge 2$ for all $m
\ge 3 m_0+4$. Thus we have the following:

\begin{thm} \label{tIn} Assume $P_{m_0}(X) \ge 2$ for some positive integer
$m_0$. If the induced map $f$ is of type $I_n$. Then
\begin{enumerate}
\item $\rho_0 \le \rho_1 \le 3m_0+4$.

\item $K^3\geq \frac{11}{12m_0(m_0+1)^2}$.

\item $\varphi_m$ is birational if $m \ge 5m_0+6$ (cf. \cite[Theorem 0.1]{JPAA}).
\end{enumerate}
\end{thm}
\end{setup}

\begin{setup}  \label{I3}{\bf Type $I_3$.}

We take $G_1=4\sigma^*(K_{S_0})$ so as to estimate $K_X^3$. Then,
as seen in \ref{In}, $\deg(K_C)\geq 18$. Being in a better
situation with $p=a_{m_0}-1\geq 2$, a better number $\beta$ can be
found. In fact, by \cite[Lemma 3.3]{Chen-Zuo}, one may take a
number sequence $\{\beta_n\}$ with $\beta_n\mapsto
\frac{p}{4(m_0+p)}\geq \frac{1}{2(m_0+2)}$ such that
$\pi^*(K_X)|_S-\beta_n C$ is numerically equivalent to an effective
$\bQ$-divisor. Namely, one may take a number $\beta\geq
\frac{1}{2(m_0+2)}-\delta$ for some $0 < \delta \ll 1$. Now
inequality (2.2) gives $\xi\geq
\frac{18}{1+\frac{m_0}{2}+\frac{1}{\beta}}$, i.e.
$\xi\geq\frac{36}{5(m_0+2)}$ by taking the limit. Hence inequality
(2.1) implies $K^3\geq \frac{36}{5m_0(m_0+2)^2}$.

We take a different $|G|$ on $S$ to study the birationality. In
fact, we will take $|G|$ to be the movable part of
$|2\sigma^*(K_{S_0})|$. A different point from previous ones is
that $|G|$ is not always base point free. But since we have the
induced fibration $f:X'\longrightarrow \Gamma$, we can consider the
relative bi-canonical map of $f$, namely the rational map $\Psi:
X'\dashrightarrow {\bf P}$ over $\Gamma$. First we can blow up the
indeterminacy of $\Psi$ on $X'$. Then we can assume, in the
birational equivalence sense, that $\Psi$ is a morphism over $B$.
By further modifying $\pi$, we can even finally assume that $\pi$
dominates $\Psi$. With this assumption (or by taking a sufficiently
good $\pi$), we see that $|G|$ is base point free since $|G|$ gives
the bicanonical morphism for each general fiber $S$ of $f$.

By Proposition \ref{nonvanishing} and Lemma \ref{a(1)}, Assumptions
\ref{asum} (1) is satisfied for $m\geq
\rounddown{\frac{5m_0}{2}}+4$. Recall that we have $p=a_{m_0}\geq
2$.

\medskip

{\bf Claim A}. Assumptions \ref{asum} (2) is satisfied for $m\geq
\text{min}\{3m_0+6, \rho_0+2m_0+2\}$.
\medskip

In fact, the argument of \ref{Ip} works here. A different place is
that we have a better bound for $\beta_n$ since $p\geq 2$, but we
only have $\deg(K_C)\geq 2$. By \cite[Lemma 3.3]{Chen-Zuo}, we can
find
 a sequence of rational numbers $\{{\beta}_n\}$ with ${\beta}_n\mapsto
\frac{p}{2(m_0+p)}$ such that $\pi^*(K_X)_{|S}-
{\beta}_n(2\sigma^*(K_{S_0}))\equiv H_n$ for effective
$\bQ$-divisors $H_n$. We may assume that  $\beta \ge
\frac{1}{m_0+2}-\delta$ for some $0<\delta \ll 1$.

Now the last three terms of inequality (2.6) can be replaced by
$$\begin{array}{rrl}
&&K_S+\roundup{(m-m_0-1)\pi^*(K_X)_{|S}}\\
&\geq& K_S+\roundup{(m-m_0-1)\pi^*(K_X)_{|S}-\frac{2}{\beta_n}H_n}\\
&=&K_S+\roundup{Q_2}+4\sigma^*(K_{S_0})
\end{array}$$
where
$Q_2:=(m-m_0-1)\pi^*(K_X)_{|S}-4\sigma^*(K_{S_0})-\frac{2}{\beta_n}H_n\equiv
(m-m_0-1-\frac{2}{\beta_n})\pi^*(K_X)|_S$ is nef and big whenever
$m\geq 3m_0+6$. According to a theorem of Xiao \cite{X}, $|G|$ is
either not composed with a pencil or composed with a rational
pencil. Thus, according to \cite[Lemma 2]{T} and Remark
\ref{separate}, Assumptions \ref{asum} (2) is satisfied for $m\geq
3m_0+6$. On the other hand, we have an inclusion:
$\OO_{\Gamma}(2)\hookrightarrow f_*\omega_{X'}^{m_0}$ which
naturally gives rise to the inclusion:
$f_*\omega_{X'/\Gamma}^2\hookrightarrow f_*\omega_{X'}^{2m_0+2}$.
Now Viehweg's semi-positivity theorem \cite{VV} implies that
$f_*\omega_{X'/\Gamma}^2$ is generated by global sections. Thus
$|(2m_0+2)K_{X'}|_{|S}$ can distinguish different generic
irreducible elements of $|G|$. So Assumptions \ref{asum} (2) is
naturally satisfied for all $m\geq \rho_0+2m_0+2$. We have proved
Claim A.
\medskip

Finally we consider the value of $\alpha$. Recall that we may take
$\beta\mapsto \frac{p}{2m_0+2p}\geq \frac{1}{m_0+2}$. Inequality
(2.2) gives $\xi\geq
\frac{2}{1+\frac{m_0}{2}+m_0+2}=\frac{4}{3(m_0+2)}$. If we take
$m=3m_0+4$. Then $\alpha>1$. Theorem \ref{technical} says $\xi\geq
\frac{4}{3m_0+4}$. Eventually, take $m\geq 3m_0+6$. Then
$\alpha>2$. Theorem \ref{technical} implies that $\varphi_m$ is
birational for all $m\geq 3m_0+6$.

We thus conclude the following:
\begin{thm} \label{tI3} Assume $P_{m_0}(X) \ge 3$ for some positive integer
$m_0$. If the induced map $f$ is of type $I_3$. Then
\begin{enumerate}
\item $\rho_0 \le \rho_1 \le \rounddown{\frac{3m_0}{2}}+4$.

\item $K^3\geq \frac{36}{5m_0(m_0+2)^2}$.

\item $\varphi_m$ is birational if $m \ge 3m_0+6$.
\end{enumerate}
\end{thm}
\end{setup}

By collecting all above results, we have the following:

\begin{cor} \label{volume} Assume $P_{m_0}(X) \ge 2$ for some positive integer
$m_0$. Then $K^3\geq \frac{11}{12m_0(m_0+1)^2}$.
\end{cor}

\begin{setup}{\bf Volume optimization.}

Indeed, when $m_0$ is small, the estimation of $K_X^3$ could be
optimized by recursively applying Theorem \ref{technical} with a
suitable $m$.

For example, suppose $m_0=11$ and $f$ is of type $III$. Then
inequality (2.2) gives $\xi\geq \frac{6}{23}$. Take $m=27$. By
Theorem \ref{technical}, we get $\xi\geq \frac{8}{27}$. So
inequality (2.1) gives $K^3\geq
\frac{8}{3267}>\frac{10}{m_0^2(3m_0+2)}$.

Let's consider another example with $m_0=8$ and $f$ being of type
$II$. Then we may take $\beta=\frac{1}{8}$. Inequality (2.2) gives
$\xi\geq \frac{2}{17}$. Take $m=26$. Then $\alpha\geq
\frac{18}{17}>1$. Theorem \ref{technical} gives $\xi\geq
\frac{2}{13}$. Take $m=24$. Then $\alpha>1$. Again, one gets
$\xi\geq \frac{1}{6}$. So inequality (2.1) implies $K^3\geq
\frac{1}{384}>\frac{4}{m_0^2(3m_0+2)}$.

With the idea mentioned above, a patient reader should have no
difficulty to check the following table on the lower bound of $K^3$
for small $m_0$. \smallskip

\centerline{\bf Table A}
\smallskip

{\normalsize
\begin{tabular}{|c|c|c|c|c|c|c|}
\hline
 $m_0$ & 2 & 3 & 4 & 5 &6 & 7\\
\hline
    $ III$ &${1}/{3}$ &${8}/{81}$ &${1}/{22}$ &${8}/{325}$ & ${1}/{72}$ &
     ${4}/{441}$
     \\ \hline
     $II$ &${1}/{8}$ &${2}/{45}$ &${1}/{52}$ &${1}/{100}$ & ${1}/{162}$ & ${4}/{1029}$\\ \hline
     $P_{m_0}\geq 3$ &${1}/{8}$ &${2}/{45}$ &${1}/{52}$ &${1}/{100}$ & ${1}/{162}$ & ${4}/{1029}$
     \\ \hline
    $P_{m_0}\geq 2$ &${5}/{96}$ &${5}/{264}$ &${1}/{108}$ &${1}/{192}$ & ${5}/{1554}$ & ${5}/{2408}$
     \\
\hline

\hline
 $m_0$& 8 & 9 & 10 & 11 &12&\\
\hline
     $III$ &
     ${1}/{160}$ & ${4}/{891}$ & ${2}/{625}$ & ${8}/{3267}$ &
     ${1}/{522}$&
     \\ \hline
     $II$
     &${1}/{384}$ & ${2}/{1053}$ & ${1}/{725}$ & ${1}/{968}$ &
     ${1}/{1224}$&
     \\ \hline
     $P_{m_0}\geq 3$
     &${1}/{384}$ & ${2}/{1053}$ & ${1}/{725}$ & ${1}/{968}$ &
     ${1}/{1224}$&
     \\ \hline
    $P_{m_0}\geq 2$
     &${5}/{3456}$ & ${1}/{954}$ & ${1}/{1276}$ & ${5}/{8448}$ & ${5}/{10764}$&    \\
\hline
\end{tabular}
}
\end{setup}

\begin{lem}\label{chi(O)>1} If $f$ is of type $I_n$ and $q(X)=0$,
then $\chi(\OO_X) \leq 1$.
\end{lem}
\begin{proof}
We have an induced fibration $f:X'\longrightarrow \Gamma$ onto the
rational curve $\Gamma$. A general fiber $S$ of $f$ is a
nonsingular projective surface of general type with $p_g(S)=0$.
Because $\chi(\OO_S)>0$, we see $q(S)=0$. This means
$f_*\omega_{X'}=0$ and $R^1f_*\omega_{X'}=0$ since they are both
torsion free by \cite{Kol}. Thus we get by \ref{inv} the following
formulae:
$$h^2(\OO_X)=h^2(\OO_{X'})=h^1(f_*\omega_{X'})+h^0(R^1f_*\omega_{X'})=0;$$
$$q(X)=q(X')=g(\Gamma)+h^1(R^1f_*\omega_{X'})=0.$$
So we see $\chi(\OO_X)=1-q(X)+h^2(\OO_X)-p_g(X)\leq 1$.
\end{proof}



\begin{setup}{\bf Miyaoka-Reid inequality on $\mathscr{B}(X)$.} We
refer to \cite[Section 2]{Explicit_I} for the definition of
baskets. Assume that Reid's basket of singularities on $X$ is
$B_X:=\mathscr{B}(X)=\{(b_i,r_i)\}$. According to \cite[10.3]{YPG},
one has
$$\frac{1}{12}K_X\cdot c_2(X)=-2\chi({\mathcal
O}_X)+\sum_i\frac{r_i^2-1}{12r_i}$$ where $c_2(X)$ is defined via
the intersection theory by taking a resolution of singularities of
$X$. On the other hand, \cite[Corollary 6.7]{Miyaoka} says
$K_X\cdot c_2(X)\geq 0$. Thus one has the following inequality
$$\sum_i r_i- 24\chi({\mathcal
O}_X)\geq \sum_i \frac{1}{r_i}. \eqno{(2.7)}$$
\end{setup}

A direct application of inequality (2.7) is the following:

\begin{cor}
Suppose that we have a packing between formal baskets: ${\bf
B}:=(B,\chi(\OO_X), \tilde{P}_2)\succcurlyeq {\bf
B'}:=(B',\chi(\OO_X), \tilde{P}_2)$ and that inequality (2.7) fails
for ${\bf B'}$. Then (2.7) fails for ${\bf B}$.
\end{cor}

\section{\bf General type 3-folds with $\chi=1$}

In this section, we always assume $\chi(\OO_X)=1$. If there is a
small number $m_0$ such that $P_{m_0}>1$, then one can detect the
birational geometry of $X$ by studying $\varphi_{m_0}$. Thus a
natural question is what a practical number $m_0$ can be found such
that $P_{m_0}>1$. This is exactly the motivation of this section.
Equivalently, we shall give a complete classification of baskets to those $X$
with $P_m\leq 1$ for $m\leq 6$.

\begin{setup}\label{p6}{\bf Assumption:} $P_m(X) \leq 1$ for $1\leq m \leq
6$.
\end{setup}

In fact, $P_m$ satisfies the following geometric condition.

\begin{lem}\label{p2m} Assume
$\chi(\mathcal{O}_X) =1$. Then $P_{m+2}\geq P_m+P_2$ for all $m\geq
2$.
\end{lem}
\begin{proof} By Reid's formula (\cite{YPG}), we have
$$P_{m+2}-P_m-P_2=(m^2+m)K_X^3-\chi(\OO_X)+(l(m+2)-l(m)-l(2)).$$
By \cite[Lemma 3.1]{Flt}, one sees $l(m+2)-l(m)-l(2)\geq 0$. Since
$K_X^3>0$ and $\chi(\OO_X)=1$, we have $P_{m+2}- P_m-P_2
> -1$.
\end{proof}

We consider the formal basket
$${\bf B}:=(B, \chi(\OO_X), P_2(X))$$
where $B=\mathscr{B}(X)$. As we have seen in \cite[Section
3]{Explicit_I},
\begin{itemize}
\item[(i)] $K^3({\bf B})=K^3(B)=K_X^3>0$;

\item[(ii)] $P_{m}({\bf B})=P_m(X)$ for all $m\geq 2$.
\end{itemize}

By Lemma \ref{p2m}, we see $P_{4}\geq 2$ if $P_2>0$. Thus under
Assumption \ref{p6}, we have $P_2=0$. We can also get
 $P_{m+2}>0$ whenever
$P_{m}>0$. Thus, in practice, we only need to study the following
types: $P_2=0$ and
$$\begin{array}{lcl} (P_3,P_4,P_5,P_6)&=&
(0,0,0,0),(0,0,0,1),(0,0,1,0),
(0,0,1,1),\\
&&(0,1,0,1),(0,1,1,1),(1,0,1,1),(1,1,1,1).
\end{array} \eqno{(3.1)}$$

Now we consider formal basket ${\bf B}:=(B, 1, 0)$. We might abuse the notation of baskets and formal baskets in this section for we always
have $\chi=1,P_2=0$ in this section. We keep the notation as in \cite{Explicit_I}.

With explicit value of $(P_3,
P_4, P_5,P_6)$, we are able to determine $\mathscr{B}^{(5)}(B)$ (cf
\cite[Sections 3,4]{Explicit_I}). Our main task is to search all
possible minimal (with regard to $\succ$) positive baskets $B_{\text
min}$ dominated by $\mathscr{B}^{(5)}(B)$. Take ${\bf
B^{5}}:=(\mathscr{B}^{(5)}(B),1,0)$ and ${\bf B}_{\text
min}:=(B_{\text min},1,0)$. Then we see: ${\bf B^{5}}\succcurlyeq
{\bf B}\succcurlyeq {\bf B}_{\text min}$.

%
%

Now we classify all minimal positive geometric baskets $B_{\text
min}$.

\begin{setup}\label{0000}{\bf Case  I}: $P_3=P_4=P_5=P_6=0$
(impossible)

We have $\sigma=10, \tau=4, \Delta^3=5, \Delta^4=14, \epsilon=0,
\sigma_5=0$ and $\epsilon_5=2$. The only possible initial basket is
 $\{5 \times (1,2), 4 \times (1,3), (1,4)\}$.
And $B^{(5)}=\{3\times (1,2), 2 \times (2,5), 2 \times (1,3),
(1,4)\}$ with $K^3=\frac{1}{60}$. We shall calculate $B_{\text min}$
of $B^{(5)}$.


If we pack $\{ (1,2),(2,5) \} $ into $\{(3,7)\}$. Then we get\\
{\bf I-1.} $B_{1,1}=\{ 2\times (1,2), (3,7),(2,5),  2\times(1,3), (1,4)  \}$, $K^3=\frac{1}{420}$   \\
which admits no further prime packing into positive baskets. Hence
$B_{1,1}$ is minimal positive.

We consider those baskets with $(1,2)$ unpacked because otherwise
it's dominated by $B_{1,1}$. So we consider the  packing:
$$\{3\times (1,2),(2,5),(3,8),(1,3),(1,4)\}$$ with $K^3=\frac{1}{120}$. This
basket
allows two further packings to minimal positive ones:\\
{\bf I-2.} $B_{1,2}=\{3\times (1,2), (2,5), (4,11) , (1,4) \}$ , $K^3=\frac{1}{220}$. \\
{\bf I-3.} $B_{1,3}=\{3\times (1,2), (5,13), (1,3) , (1,4) \}$,
$K^3=\frac{1}{156}$.

{}Finally we consider the case that both $(1,2)$ and $(2,5)$ remain
unpacked. We get one more basket which is
indeed minimal positive:\\
{\bf I-4.} $B_{1,4}=\{ 3\times (1,2),  2\times (2,5),(1,3), (2,7)
\}$ , $K^3=\frac{1}{210}$.

A direct calculation shows that $B_{1,1}$, $B_{1,2}$, $B_{1,3}$ and
$B_{1,4}$ all do not satisfy inequality (2.7). Hence $B$ does not
satisfy (2.7), a contradiction. This means that Case I is
impossible.
\end{setup}

\begin{setup}\label{0001}{\bf Case  II}: $P_3=P_4=P_5=0,P_6=1$
($\Rrightarrow B_{2,1}, B_{2,2}$)

Now we have $\sigma=10$, $\tau=4$, $\Delta^3=5$, $\Delta^4=14$,
$\epsilon \le 1$. If $\epsilon=0$, then $\epsilon_5=1$ and if
$\epsilon=1$, then $\epsilon_5=0$. Thus all possible initial baskets
and $B^{(5)}$ are as follows:

 {\bf II-i.} $B^{(0)}=\{5 \times (1,2), 4
\times (1,3), (1,4)\} \succ B^{(5)}=\{4 \times (1,2),(2,5), 3 \times
(1,3),
 (1,4)\}$, with $K^3(B^{(5)})=\frac{1}{20}$.

{\bf II-ii.}  $B^{(0)}=\{5 \times (1,2), 4 \times (1,3), (1,5)\}
\succ B^{(5)}=\{5 \times (1,2), 4 \times (1,3), (1,5)\}$,with
$K^3(B^{(5)})=\frac{1}{30}$.

In Case II-i, we first consider the situation that all single
baskets $(1,2)$ are packed into:
 $\{(6,13), 3 \times (1,3),(1,4)  \}$, which gives
 a unique minimal positive basket:\\
{\bf II-1.} $B_{2,1}=\{(6,13), (1,3),  (3,10)  \}$ , $K^3=\frac{1}{390}$, $P_9=2$, $P_{13}=3$. \\

We then consider the situation that at least one basket $(1,2)$ remains unpacked. Then we get the minimal positive basket:\\
{\bf II-2.} $B_{2,2}=\{(1,2), (5,11),  (4,13)  \}$ ,
$K^3=\frac{1}{286}$, $P_9=2$, $P_{13}=3$.

Notice, however, that if $\{3 \times (1,2),(3,7), 3 \times (1,3),
(1,4)\} \succ B$, then $B$ dominates $B_{2,2}$. Thus it remains to
consider the situation that all single baskets $(1,2)$ are unpacked,
but $(2,5)$ must be packed with some $(1,3)$. So we
get the following minimal positive baskets:\\
{\bf II-3.} $B_{2,3}=\{  (4,8), (3,8),  (3,10)  \}$ , $K^3=\frac{1}{40}$. \\
{\bf II-4.} $B_{2,4}=\{  (4,8), (4,11), (2,7) \}$ , $K^3=\frac{2}{77}$.  \\
{\bf II-5.}  $B_{2,5}=\{ (4,8), (5,14), (1,4) \}$ , $K^3=\frac{1}{28}$. \\

In Case II-ii, $B^{(5)}$ admits no further prime packing.
Thus we get:\\
{\bf II-6.} $B_{2,6}=\{ (5,10),(4,12),(1,5)\}$, $K^3=\frac{1}{30}$.

One may check that $B_{2,3}, B_{2,4}, B_{2,5}, B_{2,6}$ do not
satisfy inequality (2.7). Thus only {\bf II-1} and {\bf II-2} can
happen.
\end{setup}

\begin{setup}\label{0010}{\bf Case  III}: $P_3=P_4=0,P_5=1,P_6=0$ ($\Rrightarrow B_{3,1}\sim B_{3,5}$)

Now we have $\sigma=10$, $\tau=4$, $\Delta^3=5$, $\Delta^4=15$.
Moreover, $P_7 \ge 1$, hence $\epsilon=0$, $\sigma_5=0$ and
$\epsilon_5=4$. Thus the only possible initial basket and $B^{(5)}$
are:
 $$B^{(0)}=\{5 \times (1,2), 5 \times (1,3) \} \succ
B^{(5)}=\{(1,2),4 \times (2,5), (1,3)\}.$$
So we get minimal positive baskets:\\
{\bf III-1.} $B_{3,1}=\{(9,22),(1,3)\}$, $K^3= \frac{1}{66}$, $P_9=2$, $P_{10}=3$. \\
{\bf III-2.} $B_{3,2}=\{(7,17), (3,8)\}$, $K^3= \frac{1}{136}$, $P_{10}=2$, $P_{12}=3$. \\
{\bf III-3.} $B_{3,3}=\{(5,12), (5,13)\}$, $K^3= \frac{1}{156}$, $P_{10}=2$, $P_{12}=3$. \\
{\bf III-4.} $B_{3,4}=\{(3,7),(7,18) \}$, $K^3= \frac{1}{126}$, $P_{10}=2$, $P_{12}=3$. \\
{\bf III-5.} $B_{3,5}=\{(1,2),(9,23),\}$, $K^3= \frac{1}{46}$, $P_{8}=2$, $P_{10}=4$. \\
\end{setup}

\begin{setup}\label{0011}{\bf Case  IV}: $P_3=P_4=0,P_5=1,P_6=1$
($\Rrightarrow B_{3,1}, B_{3,2}, B_{3,4}, B_{3,5}$)

Now we have $\sigma=10$, $\tau=4$, $\Delta^3=5$, $\Delta^4=15$.
Moreover, the initial basket must have $n^0_{1,2}=n^0_{1,3}=5$,
hence $n^0_{1,r}=0$ for all $r \ge 4$. It follows that $\epsilon=0$,
$\sigma_5=0$ and $\epsilon_5=3$. Thus the only possible initial
basket and $B^{(5)}$ are:
$$B^{(0)}=\{5 \times
(1,2), 5 \times (1,3) \} \succ B^{(5)}=\{2 \times (1,2),3 \times
(2,5), 2 \times (1,3)\}.$$
So we get minimal positive baskets:\\
{\bf IV-1.} $\{(8,19),(2,6)\} \succ B_{3,1}$. \\
{\bf IV-2.} $\{(6,14), (4,11)\} \succ B_{3,4}$. \\
{\bf IV-3.} $\{(4,9), (6,16)\} \succ B_{3,2}$. \\
{\bf IV-4.} $\{(2,4),(8,21) \} \succ B_{3,5}$. \\
\end{setup}

\begin{setup}
\label{0101}{\bf Case V}: $P_3=0,P_4=1,P_5=0,P_6=1$. ($\Rrightarrow
B_{5,1}\sim B_{5,3}$)

We have $\sigma=10$, $\tau=4$, $\Delta^3=6$, $\Delta^4=13$ and
$\sigma_5 \le \epsilon \le 2$. The initial baskets have 4 types:\\
{\bf V-i.} $\{6 \times (1,2),  (1,3), 3 \times (1,4)\}$;\\
{\bf V-ii.} $\{6 \times (1,2),  (1,3), 2 \times (1,4), (1,5)\}$;\\
{\bf V-iii.} $\{6 \times (1,2),  (1,3),  (1,4),2 \times (1,5)\}$;\\
{\bf V-iv.} $\{6 \times (1,2),  (1,3), 2 \times (1,4), (1,r)\}$
with $r \ge 6$.

Cases V-iii and V-iv are impossible since $K^3 \le 0$. For Case
V-i, we have $\epsilon_5=1$ and for Case V-ii, we have
$\epsilon_5=0$. Hence $B^{(5)}$ have two possibilities,
correspondingly:\\
{\bf V-i.} $\{ (5,10), (2,5), (3,12)\}$;\\
{\bf V-ii.} $\{ (6,12),  (1,3),  (2,8), (1,5)\}$.

By computation, we get minimal positive baskets as follows:\\
{\bf V-1.} $B_{5,1}=\{(7,15),(3,12)\}$, $K^3= \frac{1}{60}$, $P_7=2$, $P_8=3$. \\
{\bf V-2.} $B_{5,2}=\{(6,12),(1,3),(3,13)\}$, $K^3= \frac{1}{39}$, $P_8=3$. \\
{\bf V-3.} $B_{5,3}=\{(6,12),(3,11),(1,5)\}$, $K^3= \frac{1}{55}$, $P_8=2$, $P_{10}=4$. \\
\end{setup}

\begin{setup}
\label{0111}{\bf Case  VI}: $P_3=0,P_4=P_5=P_6=1$ ($\Rrightarrow
B_{6,1}\sim B_{6,6}$)

We have $\sigma=10$, $\tau=4$, $\Delta^3=6$, $\Delta^4=14$. Also
$P_7 \ge 1$ and hence $\sigma_5 \le \epsilon \le 2$. The initial
baskets have 4 types:\\
{\bf VI-i.} $\{6 \times (1,2),  2 \times (1,3), 2 \times
(1,4)\}$;\\
{\bf VI-ii.} $\{6 \times (1,2),  2 \times (1,3), (1,4), (1,5)\}$;\\
{\bf VI-iii.} $\{6 \times (1,2), 2 \times  (1,3),  2 \times
(1,5)\}$;\\
{\bf VI-iv.} $\{6 \times (1,2), 2 \times  (1,3), (1,4), (1,r)\}$
with $r \ge 6$.

Since there are only 2 baskets of $(1,3)$, we have $\epsilon_5
=3-\sigma_5 \le 2$. Hence $\sigma_5 >0$ and $\epsilon >0$.
Therefore, Case VI-i is impossible.

For Case VI-ii, $\epsilon_5=2$, hence\\
{\bf VI-ii.} $B^{(5)}=\{ 4 \times (1,2),2 \times (2,5), (1,4),(1,5)
 \}$.

We get minimal positive baskets as follows:\\
{\bf VI-1.} $B_{6,1}=\{ (1,2),(7,16), (2,9)\}$, $K^3=\frac{1}{144}$, $P_7=2$, $P_9=3$.\\
{\bf VI-2.} $B_{6,2}=\{ (6,13),(2,5), (2,9)\}$, $K^3=\frac{8}{585}$, $P_7=2$, $P_8=3$.\\
{\bf VI-3.} $B_{6,3}=\{ (8,18),(1,4), (1,5)\}$, $K^3=\frac{1}{180}$, $P_7=2$, $P_9=3$.\\

For Case VI-iii, $\epsilon_5=1$, hence\\
{\bf VI-ii.} $B^{(5)}=\{ 5 \times (1,2), (2,5),(1,3) , 2 \times
(1,5)
 \}$.

Then we get minimal positive baskets as follows:\\
{\bf VI-4.} $B_{6,4}=\{ (1,2),(6,13),(1,3),(2,10)\}$, $K^3=\frac{1}{390}$, $P_8=2$, $P_9=3$.\\
{\bf VI-5.} $B_{6,5}=\{ (5,10),(3,8), (2,10)\}$, $K^3=\frac{1}{40}$, $P_8=3$.\\

For Case VI-iv, $\epsilon_5=2$, hence\\
{\bf VI-iv.} $B^{(5)}=\{ 4 \times (1,2),2 \times (2,5), (1,4),
(1,r) \}$ with $r \ge 6$.

 Since $K^3(B^{(5)})>0$, we must have
 $r=6$. Then we get the minimal positive basket:\\
{\bf VI-6.} $B_{6,6}=\{ (3,6),(3,7),(2,5), (1,4),(1,6)\}$, $K^3=\frac{1}{420}$, $P_{10}=2$, $P_{12}=3$.\\
\end{setup}

\begin{setup}
\label{1011}{\bf Case  VII}: $P_3=1,P_4=0,P_5=P_6=1$ (impossible)

We have $\sigma=9$, $\tau=3$, $\Delta^3=1$, $\Delta^4=9$. Moreover,
$P_7 \geq 1$ and hence $\epsilon=0$. It follows that $\sigma_5=0$
and $\epsilon_5=2$. The initial basket is: $B^{(0)}=\{ (1,2), 7
\times (1,3),(1,4) \}.$

 Note that there is only one basket of type $(1,2)$. However,
 Since $\epsilon_5=2$, one has $1\geq n_{2,5}^5=2$, a contradiction.
 Thus Case VII does not happen.
\end{setup}

\begin{setup}
\label{1111}{\bf Case  VIII}: $P_3=P_4=P_5=P_6=1$ ($\Rrightarrow
B_{8,1}\sim B_{8,3}$)

We have $\sigma=9$, $\tau=3$, $\Delta^3=2$, $\Delta^4=8$. Moreover,
$P_7 \geq 1$ and then $\epsilon \le 1$. If $\epsilon=1$, then
$\sigma_5=1$ and $\epsilon_5=1$. If $\epsilon=0$, then $\sigma_5=0$
and $\epsilon_5=2$. The initial baskets and $B^{(5)}$ have 2
types:\\
{\bf VIII-i.} $B^{(0)}=\{2 \times (1,2), 4 \times (1,3), 3 \times
(1,4)\} \succ B^{(5)}=\{2 \times (2,5), 2 \times (1,3), 3 \times
 (1,4)\}$ with $K^3(B^{(5)})=\frac{1}{60}$.\\
{\bf VIII-ii.}  $B^{(0)}=\{2 \times (1,2), 4 \times (1,3), 2 \times
(1,4),(1,5)\} \succ B^{(5)}=\{(1,2)$, $(2,5), 3 \times (1,3), 2
\times (1,4), (1,5)\}$ with $K^3(B^{(5)})=0$.

Clearly, Case VIII-ii is impossible since $K^3$ is not positive.

For Case VIII-i, we first consider the situation that one single
basket $(2,5)$ is packed, so that we get  the basket:
$\{(2,5),(3,8),(1,3),(4,12)\}$. We can get two minimal positive
baskets as follows:\\
{\bf VIII-1.} $B_{8,1}= \{(5,13),(1,3),(3,12)\}$, $K^3=\frac{1}{156}$, $P_7=2$, $P_8=3$.\\
{\bf VIII-2.} $B_{8,2}= \{(2,5),(4,11),(3,12)\}$, $K^3=\frac{1}{220}$, $P_7=2$, $P_8=3$.\\

It remains to consider the situation that each single basket $(2,5)$
remains unpacked. We then obtain the basket:
$$B_{210}:=\{(4,10),(1,3),(2,7),(2,8)\}$$ with $K^3=\frac{1}{210}$, $P_7=2$, $P_{10}=3$.
After a one-step prime packing, we get the minimal positive basket:\\
{\bf VIII-3.} $B_{8,3}=\{(4,10),(1,3),(3,11),(1,4)\}$,
$K^3=\frac{1}{660}$, $P_7=2$.
\end{setup}

The detailed classification (3.3$\sim$ 3.10) makes it possible for
us to study the birational geometry of $X$, of which the first
application is the following theorem.

\begin{thm}\label{420}\footnote{With a different approach,
L. Zhu \cite{Zhu} also proved $K^3\geq \frac{1}{420}$.} Assume
$\chi(\OO_X)=1$. Then $K_X^3 \geq \frac{1}{420}$. Furthermore,
$K_X^3=\frac{1}{420}$ if, and only if, $\mathscr{B}=B_{6,6}$.
\end{thm}
\begin{proof} If $\mu_1\leq 6$, then Proposition \ref{volume}
implies $K_X^3\geq \frac{1}{294}\cdot \frac{11}{12}>\frac{1}{420}$.

We may assume that $P_m \leq 1 $ for $m \le 6$. We have seen
$P_2=0$. Since ${\bf B^{5}}\succcurlyeq {\bf B}\succcurlyeq {\bf
B}_{\text min}$ and by \cite[Lemma 3.6]{Explicit_I}, we have
$$K_X^3=K^3(B)\geq K^3(B_{\text{min}})$$
where $B_{\text{min}}$ is in the set $\{B_{2,1}, B_{2,2},
B_{3,1}\sim B_{3,5}, B_{5,1}\sim B_{5,3}, B_{6,1}\sim B_{6,6},
B_{8,1}\sim B_{8,3}\}.$

If $B_{\text{min}}\neq B_{6,6}, B_{8,3}$, then we have seen
$K^3(B_{\text{min}})>\frac{1}{420}$.

If $B_{\text{min}}=B_{8,3}$, we show $B\neq B_{8,3}$. In fact,
if $B=B_{8,3}$, then $P_7(B)=2$ as we have seen in 3.10. By Table A
in Section 2, we have $K_X^3=K^3({\bf B}) \geq \frac{5}{2408} >
\frac{1}{660}$, a contradiction. Hence $B\succ B_{8,3}$. Notice
that $B_{8,3}$ is obtained, exactly, by one-step packing from
$$B_{210}:=\{(4,10),(1,3),(2,7),(2,8)\}$$ and no other ways.
This says $B\succcurlyeq B_{210}$ and so $K_X^3\geq
K^3(B_{210})=\frac{1}{210}$.

We have seen $K^3(B_{6,6})=\frac{1}{420}$. We are done.
\end{proof}

The proof of the last theorem gives the following:

\begin{cor}\label{min} Assume $\chi(\OO_X)=1$ and $P_m\leq 1$ for all $m\leq 6$. Then
$\mathscr{B}(X)$ either dominates a minimal basket in the set
$$\{B_{2,1}, B_{2,2}, B_{3,1}\sim B_{3,5}, B_{5,1}\sim B_{5,3},
B_{6,1}\sim B_{6,6}, B_{8,1}, B_{8,2}\}$$ or dominates the basket
$B_{210}$.
\end{cor}


\begin{cor}\label{10} Assume $\chi(\OO_X)=1$. Then $P_{10}(X)\geq
2$ and, in particular, $\mu_1\leq 10$.
\end{cor}
\begin{proof} If $P_{m_0}\geq 2$ for some $m_0\leq 6$, then, by
Lemma \ref{p2m}, one can see $P_{10}\geq 2$. Otherwise, Corollary
\ref{min} and \cite[Lemma 3.6]{Explicit_I} imply that
$P_{10}=P_{10}(\mathscr{B}(X))\geq P_{10}(B_{*})$ where $B_{*}$
denotes a minimal positive basket mentioned in Corollary \ref{min}.
By a direct computation, we get $P_{10}(B_{*})\geq 2$.
\end{proof}

Example \ref{ex1} shows that the statement in Corollary \ref{10} is
optimal since $P_9(X_{46})=1$.

\begin{thm}\label{7}\footnote{D. Shin \cite{Shin} proved the first statement
in a different way} Assume $\chi(\OO_X)=1$. Then
\begin{itemize}
\item[(1)]
$\rho_0\leq  7$.

\item[(2)] Either $P_5>0$ or $P_6>0$.
\end{itemize}
\end{thm}
\begin{proof} (1) Recall that $\mu_0:=\text{min}\{m|P_m>0\}$. By 3.3,
we see $\mu_0\leq 6$.

When $\mu_0\leq 3$, it is easy to deduce the statement by Lemma
\ref{p2m}.

When $\mu_0=4$, Lemma \ref{p2m} implies $P_{2k}>0$ for all $k\geq 3$.
If $P_7>0$, Lemma \ref{p2m} implies $P_{2k+1}>0$ for all $k\geq 3$
and the statement (1) is true. Assume $P_7=0$. Then $P_5=0$. Now
$\epsilon_5=2-P_6-\sigma_5\geq 0$ implies $\sigma_5\leq 2-P_6\leq
1$. On the other hand, $\epsilon_6=P_4+P_6-\epsilon=0$ implies
$\epsilon\geq 2$. This means $\sigma_5=P_6=P_4=1$ and the situation
corresponds to 3.7. Thus $B\succcurlyeq B_{\text{min}}$ where
$B_{\text{min}}=B_{5,2}, B_{5,3}$. But the computation tells
$P_7(B_{\text{min}})>0$, a contradiction.

When $\mu_0=5$, we study $P_8$. If $P_8>0$, then (1) is true by Lemma
\ref{p2m}. Assume $P_8=0$. Then $P_6=0$. Now
$\epsilon_6=P_5-P_7-\epsilon=0$ gives $\epsilon=0$ and $P_5=P_7$
since $P_7\geq P_5$. Since $n_{1,4}^0=1-P_5\geq 0$, we see $P_5=1$.
So the situation corresponds to 3.5. Since the computation shows
$P_8\geq P_8(B_{3,*})>0$, a contradiction.

Finally, when $\mu_0=6$, we study $P_7$. If $P_7>0$, then Lemma
\ref{p2m} implies (1). Otherwise, $P_7=0$. Now
$\epsilon_6=P_6-\epsilon=0$ implies $\epsilon=P_6>0$. Besides,
$\epsilon_5=2-P_6-\sigma_5\geq 0$
says $P_6\leq 1$ since $\sigma_5>0$. Hence $\epsilon=P_6=1$. The
situation corresponds to 3.4. But the computation shows $P_7\geq
P_7(B_{2,1})>0$ or $P_7\geq P_7(B_{2,2})>0$, a contradiction.

(2) Assume $P_5=P_6=0$. Then Lemma \ref{p2m} implies $P_3=P_4=0$.
The situation corresponds to 3.3, which is impossible as already
seen there.
\end{proof}

\section{\bf General type 3-folds with $\chi>1$}


In this section, we assume $\chi(\OO_X)>1$. Again, we will
frequently apply our formulae and inequalities in \cite[Sections
3,4]{Explicit_I}.

When $P_{m_0}\geq 2$ for some positive integer $m_0\leq 12$, known
theorems will give an effective lower bound of $K_X^3$ and a
practical pluricanonical birationality. Therefore, similar to
Section 3, we need to classify $X$ up to baskets when preceding
plurigenera are smaller. For this reason, we make the following:

\begin{setup} \label{12}
{\bf Assumption:} $P_m \leq 1$ for all $m \leq 12$.
\end{setup}

According to \cite[Lemma 4.8]{Explicit_I}, we have seen that
$P_2=0$ under Assumption \ref{12}. Note that inequality
\cite[(3.14)]{Explicit_I}, for general type 3-folds, is as follows:
$$2P_5+3P_6+P_8+P_{10}+P_{12} \ge \chi +10
P_2+4P_3+P_7+P_{11}+P_{13}+R \eqno{(4.1)}$$ where { $$
\begin{array}{lll}
R&:=&14 \sigma_5-12 n^0_{1,5}-9
n^0_{1,6}-8n^0_{1,7}-6n^0_{1,8}-4n^0_{1,9}-2n^0_{1,10}-n^0_{1,11}\\
 &=&
2n^0_{1,5}+5n^0_{1,6}+6n^0_{1,7}+8n^0_{1,8}+10n^0_{1,9}+12n^0_{1,10}+13n^0_{1,11}\\
&&+14 \sum_{r \ge 12} n^0_{1,r}. \end{array}$$} and
$\sigma_5=\sum_{r \geq 5} n^0_{1,r}$.

Inequality (4.1) and Assumption \ref{12} implies that both $\chi$
and $P_{13}$ are upper bounded. Thus our formulae in \cite[Section
4]{Explicit_I} allow us to explicitly compute $B^{(12)}$. To be more
solid, we prove the following:

\begin{prop}\label{r>5} Assume $\chi(\OO_X)>1$ and $P_m\leq 1$ for
all $m\leq 12$. Then the formal basket ${\bf B}={\bf
B}(X):=(\mathscr{B}(X), \chi(\OO_X), 0)$ has a finite number of
possibilities.
\end{prop}

\begin{proof} We study $n_{1,r}^0$ for
$r\geq 6$. If there exists a number $r\geq 6$ such that $n^0_{1,r}
\ne 0$, then $R \geq 5$ by the definition of $R$ in inequality
(4.1). Hence, by (4.1), one has
$$8\geq 2P_5+3P_6+P_8+P_{10}+P_{12}  \geq \chi+5 \geq 7.$$
This implies that  $P_5=P_6=1$. Hence $P_{11}=1$. Now (4.1) again
reads: $5+P_8+P_{10}+P_{12} \ge 8+P_7+P_{13}$. It follows that
$P_8=P_{10}=P_{12}=1$ and $P_7=P_{13}=0$. This gives a
contradiction since $P_{13} \ge P_5 P_8=1$. So we conclude
$n^0_{1,r}=0$ for all $r\geq 6$. In other words, \cite[Assumption
3.8]{Explicit_I} is satisfied.

This essentially allows us to utilize those formulae in the last
part of \cite[Section 3]{Explicit_I}. In particular, one sees that
each quantity there is bounded and hence $B^{(12)}$ has a finite
number of possibilities. Dominated by $B^{(12)}$ (i.e.
$B^{(12)}\succcurlyeq B$), $B=\mathscr{B}(X)$ also has a finite
number of possibilities. We are done.
\end{proof}

\begin{setup}{\bf Complete classification of ${\bf B}$ satisfying Assumption \ref{12}.}

Note that, for all $0<m,n\leq 12$, and $m+n \le 13$,
$$P_{m+n} \geq P_m P_n \eqno{(4.2)}$$
naturally holds since $P_m,\ P_n\leq 1$.

Suppose we have known $B^{(12)}$. Then we can determine all
possible minimal positive baskets $B_{\text{min}}$ dominated by
$B^{(12)}$, where $B_{\text{min}}\in T$ (a finite set). Now the
formal basket ${\bf B}$ satisfies the following relation:
$$
(B^{(12)}, \chi, 0)\succcurlyeq {\bf B}\succcurlyeq
(B_{\text{min}}, \chi, 0)$$ for some $B_{\text{min}}\in T$.
Therefore, by \cite[Lemma 3.6]{Explicit_I}, we have
$K^3_X=K^3(B)\geq K^3(B_{\text{min}})>0$ and $P_m=P_m(B)\geq P_m(
B_{\text{min}})$. This is the whole strategy.

The calculation can be done by a simple computer program, or even
by a direct handy work. Our main result is Table C which is a
complete list of all possibilities of $B^{(12)}$ and its minimal
positive elements.

In fact, first we preset $P_m=0$, $1$ for $m=3,\cdots, 11$. Then
$\epsilon_6=0$ gives the value of $\epsilon$. So we know the value
of $n_{1,5}^0$. By inequality (4.1) we get the upper bound of
$\chi$ since $P_{13}\geq 0$. Since $n_{1,4}^7\geq 0$, we get the
upper bound of $\eta$. Similarly $n^9_{2,9}\geq 0$ gives the upper
bound of $\zeta$. Also $n_{4,9}^{11}\geq 0$ yields $\alpha\leq
\zeta$. Finally $n_{3,8}^{11}\geq 0$ gives the upper bound of
$\beta$. Now we set $P_{12}=0$, $1$. Then inequality (4.1) again
gives the upper bound of $P_{13}$, noting that $\chi\geq 2$.
Clearly there are, at most, finitely many solutions. With
inequality (4.2) imposed, we can get about 80 cases. An important
property to mention is the inequality: $K^3(B^{(12)})\geq
K^3(B)=K_X^3>0$. With $K^3>0$ imposed on, we have got 63 outputs,
which is exactly Table C. Simultaneously, we have been able to
calculate all those minimal positive baskets dominated by
${B}^{(12)}$, since ${B}^{(12)}$ is ``nearly'' minimal in most
cases.


If one would like to take a direct calculation by hands, it is of
course possible. Consider no. 2 case in Table C as an example.
Since $P_2=0$, $P_3=\cdots =P_7=0$, $P_8=1$ and $P_9=P_{10}
=P_{11}=0$, \cite[(3.10)]{Explicit_I} tells $\epsilon=0$ and thus
$\sigma_5=0$, which means $R=0$. Now inequality (4.1) gives
$P_{12}+1\geq \chi+P_{13}\geq 2$. So $P_{12}=1$, $\chi=2$ and
$P_{13}=0$. Now the formula for $\epsilon_{10}$ gives
$\epsilon_{10}=-\eta\geq 0$, which means $\eta=0$. Similarly
$n^9_{1,5}=\zeta-1\geq 0$. On the other hand,
$n_{3,7}^9=1-\zeta\geq 0$. Thus $\zeta=1$. Now
$n_{4,9}^{11}=\zeta-\alpha\geq 0$ gives $\alpha \leq 1$.
$n^{11}_{3,11}=1-\zeta-\alpha-\beta\geq$ gives $\alpha=\beta=0$.
Finally we get
$$\{n_{1,2},n_{5,12},...,n_{1,5}\}=\{4,0,1,0,0,2,1,0,3,0,0,0,2,0,0\}$$
That is $B^{(12)}=\{4 \times(1,2), (4,9),2 \times (2,5),(3,8),3
\times (1,3), 2 \times (1,4)\}$.

We see that $B^{(12)}$ admits only one prime packing of type
$$\{(2,5),(3,8)\} \succ \{(5,13)\}$$ over the minimal positive basket
$\{4 \times(1,2), (4,9), (2,5),(5,13),3 \times (1,3), 2 \times
(1,4)\}$. We simply write this as $\{(5,13),*\}$ in Table C. It is
now easy to calculate $K^3$ for both $B^{(12)}$ and the minimal
positive basket $\{(5,13),*\}$. Finally we can directly calculate
$P_{m}$. At the same time, $\mu_1$ is given in the table. For our
needs in the context, we also display the value of
$P_{18}=P_{18}(B^{(12)})$ or $P_{18}(B_{\text{min}})$ and
$P_{24}=P_{24}(B^{(12)})$ or $P_{24}(B_{\text{min}})$ in Table C,
though the symbols $P_{18}$ or $P_{24}$ are misused here.

So theoretically we can finish our classification by detailed
computations. We omit the details because all calculations are
similar.
\end{setup}

\begin{setup}{\bf Notation.} By abuse of notation, we denote by
$B_{*}$ the final basket corresponding to No.* in Table C. For
example, $B_2=\{4\times (1,2), (4,9), 2\times (2,5), (3,8), 3\times
(1,3), 2\times (1,4)\}$ while $B_{2a}=\{4\times (1,2), (4,9),
(2,5), (5,13), 3\times (1,3), 2\times (1,4)\}$ is minimal positive.
The relation is as follows:
$$B_2\succcurlyeq B\succcurlyeq B_{2a}.$$
Clearly, for this case, we have $\frac{1}{360}=K^3(B_{2})\geq
K_X^3\geq K^3(B_{2a})=\frac{1}{1170}.$

Another typical example is No.63, where we have $$B_{63}=\{5\times
(1,2), (4,9),  2\times (3,7), (2,5), (3,8), (4,11), 3\times (1,3),
(2,7), (1,5)\}$$ which is already minimal positive. So we have the
relation:
$$B^{(12)}=B_{63}=B=B_{\text{min}}$$
and thus $K_X^3=\frac{1}{5544}$. Of course, we will see that No.63
does not happen on any $X$.
\end{setup}

Now we begin to analyze Table C and pick out ``impossible'' cases.

\begin{prop}\label{impos} In Table C, $B\neq B_{*}$ for any $B_{*}$ in the set
$$\{B_{4a}, B_{9}, B_{16a}, B_{16c}, B_{18a}, B_{20a}, B_{21a},
B_{22}, B_{24}, B_{27a}, $$ $$B_{29a}, B_{33a}, B_{44b}, B_{46a},
B_{47}, B_{52a}, B_{55}, B_{60a}, B_{61}, B_{63}\}.$$ In
particular, cases No. 9, No. 22, No. 24, No. 47, No. 55, No. 61 and
No. 63 do not happen at all.
\end{prop}

\newpage

\centerline{\bf Table C}

{\scriptsize
$$
\begin{array}{lcccccccc}
No. & (P_3,\cdots,P_{11}) &P_{18}&P_{24}&\mu_1 & \chi   &
B^{(12)}=(n_{1,2},n_{5,11},\cdots,n_{1,5}) \text{ or }
B_{min} & K^3\\
\hline

1 &(0,0,0, 0, 0, 0, 0, 1, 0) &4&8&14&  2& (5,0,0,1,0,3,0,0,3,0,0,1,0,0,0) & \frac{3}{770}\\
2 &( 0, 0, 0, 0, 0, 1, 0, 0, 0)&3&7&15&  2 &( 4, 0,1, 0,0, 2, 1,0, 3, 0, 0,0, 2,  0, 0) & \frac{1}{360} \\
2a&&2&3 &18&& \{(2,5),(3,8),*\} \succ\{(5,13),* \} & \frac{1}{1170}\\
3&( 0, 0, 0, 0, 0, 1, 0, 1, 0 ) &3&7&15&  3&(6,1,0, 0,0,4, 1,0,4, 0, 1,0,2,0,0) & \frac{23}{9240} \\
3a&&2&3 &18&&\{(2,5),(3,8),*\} \succ \{(5,13),* \} & \frac{17}{30030}\\
4&( 0, 0, 0, 0, 0, 1, 0, 1, 0 ) &4&9&14&  3&(7,0,1, 0,0,4, 0,1,3, 0, 1,0,2,0,0) & \frac{13}{3465} \\
4a&&1&2 &14&&\{(4,11),(2,6),*\} \succ \{(6,17),* \} & \frac{1}{5355}\\
5&( 0, 0, 0, 0, 0, 1, 0, 1, 0 ) &5&10&14&  3&(7,0,1, 0,0,4, 1,0,4, 0, 0,1,1,0,0) & \frac{17}{3960} \\
5a&&4&3&15&&\{(8,20),(3,8),*\} \succ \{(11,28),* \} & \frac{1}{1386}\\
5b&&3&3&15&& \{(5,13),(4,15),* \} & \frac{1}{1170}\\
6& (0, 0, 0, 1, 0, 0, 0, 1, 0 ) &3&6&14& 3 & (9,0,0, 2,0,1, 0,1,4, 0,2,0,0, 0, 1) & \frac{1}{462} \\
7& (0, 0, 0, 1, 0, 0, 1, 0, 0 ) &3&5&14& 2 & (5,0,1, 1,0,0, 0,0,5, 0,1,0,0, 0, 1) & \frac{1}{630} \\
7a&&2&3&14&&\{(4,9),(3,7),*\} \succ \{(7,16),* \} & \frac{1}{1680}\\
8& (0, 0, 0, 1, 0, 0, 1, 1, 0 ) &3&5&14& 3 & (7,1,0, 1,0,2, 0,0,6, 0,2,0,0, 0, 1) & \frac{1}{770} \\
9& (0, 0, 0, 1, 0, 1, 0, 0, 0 ) &2&2&14& 3 & (9,0,0, 2,0,0, 1,1,4, 0,1,0,0, 1, 0) & \frac{1}{5544} \\
10& (0, 0, 0, 1, 0, 1, 0, 0, 0 ) &3&6&14& 3 & (8,0,1, 1,0,0, 2,0,5, 0,1,0,1, 0, 1) & \frac{1}{630} \\
10a&&2&4&14&&\{(4,9),(3,7),*\} \succ \{(7,16),* \} & \frac{1}{1680}\\
11& (0, 0, 0, 1, 0, 1, 0, 1, 0 ) &2&4&14& 3 & (9,0,0, 2,0,0, 1,1,3, 1,0,0,1, 0, 1) & \frac{3}{3080}\\
11a&&2&3&14&&\{(3,8),(4,11),*\} \succ \{(7,19),* \} & \frac{1}{2660}\\
12& (0, 0, 0, 1, 0, 1, 0, 1, 0 ) &5&11&14& 3 & (9,0,1, 0,0,1, 2,0,4, 0,2,0,0, 0, 1) & \frac{1}{252} \\
12a&&4&6&14&&\{(2,5),(6,16),*\} \succ \{(8,21),* \} & \frac{1}{630}\\
13& (0, 0, 0, 1, 0, 1, 0, 1, 0 ) &3&4&14& 4 & (12,0,0, 2,0,2, 0,2,4, 0,2,0,0,1,0) & \frac{4}{3465} \\
14& (0, 0, 0, 1, 0, 1, 0, 1, 0 ) &3&6&14& 4 & (10,1,0, 1,0,2, 2,0,6,0,2,0,1, 0, 1) & \frac{1}{770} \\
15& (0, 0, 0, 1, 0, 1, 0, 1, 0 ) &4&8&14& 4 & (11,0,1, 1,0,2, 1,1,5,0,2,0,1, 0, 1) & \frac{71}{27720} \\
15a&&2&4&14&&\{(4,11),(1,3),*\} \succ \{(5,14),* \} & \frac{1}{2520}\\
15b&&3&4&14&&\{(2,5),(3,8),*\} \succ \{(5,13),* \} & \frac{23}{36036}\\
15c&&3&5&14&& \{(7,16),(7,19),* \} & \frac{31}{31920}\\
16& (0, 0, 0, 1, 0, 1, 0, 1, 0 ) &5&9&14& 4 & (11,0,1, 1,0,2, 2,0,6,0,1,1,0, 0, 1) & \frac{43}{13860} \\
16a&&4&3&14&&\{(4,10),(3,8),*\} \succ \{(7,18),* \} & \frac{1}{3080}\\
16b&&4&4&14&&\{(2,5),(6,16),*\} \succ \{(8,21),* \} & \frac{1}{1386}\\
16c&&3&3&14&& \{(7,16),(5,13),* \} & \frac{3}{16016}\\
17& (0, 0, 0, 1, 0, 1, 0, 1, 1 ) &3&6&14& 3 & (9,0,0, 2,0,0, 0,2,3, 0, 1,0,1,0, 1) &  \frac{3}{1540} \\
18& (0, 0, 0, 1, 0, 1, 0, 1, 1 ) &4&7&14& 3 & (9,0,0, 2,0,0, 1,1,4, 0, 0,1,0,0, 1) &  \frac{23}{9240} \\
18a&&2&3&14&&\{(4,11),(1,3),*\} \succ \{(5,14),* \} & \frac{1}{3080}\\
18b&&4&6&14&&\{(3,8),(4,11),*\} \succ \{(7,19),* \} & \frac{83}{43890}\\

19& (0, 0, 0, 1, 0, 1, 1, 0, 0 ) &3&3&14& 3 & (8,0,1,1,0,1,0,1,5,0,1,0,0,1,0) & \frac{2}{3465} \\
20&(0, 0, 0, 1, 0, 1, 1, 0, 0 ) &4&7&14& 3 & (7,0,2, 0,0,1, 1,0,6, 0,1,0,1,0, 1) & \frac{1}{504} \\
20a&&3&3&18&&\{(2,5),(3,8),*\} \succ \{(5,13),* \} & \frac{1}{16380}\\
21&( 0, 0, 0, 1, 0, 1, 1, 1, 0 ) &4&8&14& 2 & (6,0,1, 0,0,0, 1,0,3, 1,0,0,0, 0, 1 ) & \frac{1}{360} \\
21a&&2&3&16&&\{(1,3),(3,10),*\} \succ \{(4,13),* \} & \frac{1}{4680}\\
22&(0, 0, 0, 1, 0, 1, 1, 1, 0 ) &2&3&18& 3 & (7,1,0, 1,0,1, 1,0,5, 1,0,0,1,0,1) & \frac{1}{9240}\\
23&(0, 0, 0, 1, 0, 1, 1, 1, 0 ) &3&5&14& 3 & (8,0,1, 1,0,1, 0,1,4, 1,0,0,1,0,1) & \frac{19}{13860}\\
23a&&2&3&14&&\{(4,9),(3,7),*\} \succ \{(7,16),* \} & \frac{1}{2640}\\
24&(0, 0, 0, 1, 0, 1, 1, 1, 0 ) &3&3&14& 4 & (10,1,0, 1,0,3, 0,1,6, 0,2,0,0,1,0)& \frac{1}{3465} \\
25&(0, 0, 0, 1, 0, 1, 1, 1, 0 ) &4&7&14& 4 & (9,1,1, 0,0,3, 1,0,7, 0,2,0,1,0,1) & \frac{47}{27720} \\
25a&&4&6&14&&\{(5,11),(4,9),*\} \succ \{(9,20),* \} & \frac{1}{840}\\
26&(0, 0, 0, 1, 0, 1, 1, 1, 0, ) &5&9&14& 4 & (10,0,2, 0,0,3, 0,1,6,0,2,0,1,0,1) & \frac{41}{13860} \\
26a&&3&5&14&&\{(4,11),(1,3),*\} \succ \{(5,14),* \} & \frac{1}{1260}\\
27&(0, 0, 0, 1, 0, 1, 1, 1, 0 ) &6&10&14& 4 & (10,0,2, 0,0,3, 1,0,7, 0,1,1,0,0,1) & \frac{97}{27720} \\
27a&&5&3&14&&\{(6,15),(3,8),*\} \succ \{(9,23),* \} & \frac{19}{79695}\\
27b&&5&5&14&& \{(5,13),(5,18),* \} & \frac{1}{1170}\\
28&(0, 0, 0, 1, 0, 1, 1, 1, 1 ) &4&8&14& 2 & (5,1,0, 0,0,0, 1,0,4, 0,1,0,0,0,1) & \frac{23}{9240}\\
29&(0, 0, 0, 1, 0, 1, 1, 1, 1 ) &5&10&14& 2 & (6,0,1, 0,0,0, 0,1,3, 0,1,0,0,0,1) & \frac{13}{3465}\\
29a&&2&3&14&&\{(4,11),(2,6),*\} \succ \{(6,17),* \} & \frac{1}{5355}\\
\end{array}
$$

\newpage
$$
\begin{array}{lccccccc}
No. & (P_3,\cdots,P_{11}) &P_{18}&P_{24}&\mu_1& \chi & (n_{1,2},n_{4,9},\cdots,n_{1,5}) \text{ or } B_{min} & K^3\\
\hline
30&(0, 0, 0, 1, 0, 1, 1, 1, 1 ) &3&5&14& 3 & (7,1,0, 1,0,1, 0,1,5, 0,1,0,1,0, 1) & \frac{1}{924} \\
31&(0, 0, 0, 1, 0, 1, 1, 1, 1 ) &4&6&14& 3 & (7,1,0, 1,0,1, 1,0,6, 0,0,1,0,0, 1) & \frac{1}{616} \\
32&(0, 0, 0, 1, 0, 1, 1, 1, 1 ) &5&8&14& 3 & (8,0,1, 1,0,1, 0,1,5, 0,0,1,0,0, 1) & \frac{2}{693} \\
32a&&4&6&14&&\{(4,9),(3,7),*\} \succ \{(7,16),* \} & \frac{1}{528}\\
32b&&2&2&14&&\{(4,11),(1,3),*\} \succ \{(5,14),* \} & \frac{1}{1386}\\
33 & (  0, 0, 0, 1, 1, 0, 0, 1, 0 ) &2&4&14 & 2 & (5,0, 0, 2,0,0,1,0,1,1,1,0, 0, 0, 0) & \frac{1}{840} \\
33a&&1&3&14&& \{(3,10),(2,7),*\} \succ \{(5,17),* \} & \frac{1}{2856}\\
34 & ( 0, 0, 0, 1, 1, 0, 0, 1, 0 ) &4&8& 14& 3 &( 7,0, 1, 1,0, 2, 1,0,3,0, 3,0,0,  0, 0) & \frac{1}{360} \\
34a&&3&6 &14&&\{(4,9),(3,7),*\} \succ \{(7,16),* \} & \frac{1}{560}\\
34b&&3& 4&14&&\{(2,5),(3,8),*\} \succ \{(5,13),* \} & \frac{1}{1170}\\
%
35& (  0, 0, 0, 1, 1, 0, 0, 1, 1 ) &3&6& 14& 2 &( 5,0, 0, 2,0,0,0,1,1,0,2,0,0, 0, 0) & \frac{1}{462}\\
36 & ( 0, 0, 0, 1, 1, 0, 1, 1, 0 ) &3&5&14& 2 & (4,0,1, 1,0,1, 0,0,2,1,1,0,0,  0, 0) & \frac{1}{630}\\
36a&&2&3&14&& \{(4,9),(3,7),*\} \succ\{(7,16),* \} & \frac{1}{1680}\\
36b&&2&4&14&&  \{(3,10),(2,7),*\} \succ\{(5,17),* \} & \frac{4}{5355}\\

37 & ( 0, 0, 0, 1, 1, 0, 1, 1, 0 ) &5&9&14 & 3 & (6,0,2, 0,0,3, 0,0,4, 0,3,0,0, 0, 0)&  \frac{1}{315} \\

38 & ( 0, 0, 0, 1, 1, 0, 1, 1, 1 ) &3&5&14& 2 &( 3,1,0, 1,0,1, 0,0,3,0,2,0,0, 0, 0) & \frac{1}{770}\\
39 & ( 0, 0, 0, 1, 1, 1, 0, 1, 0) &3&6&14 & 3 & (7,0,1, 1,0,1, 2,0,2,1,1,0,1,  0, 0) & \frac{1}{630} \\
39a&&2&4 &14&&\{(4,9),(3,7),*\} \succ \{(7,16),* \} & \frac{1}{1680}\\
39b&&2&5 &14&& \{(3,10),(2,7),*\} \succ \{(5,17),* \} & \frac{4}{5355}\\
40 & ( 0, 0, 0, 1, 1, 1, 0, 1, 0) &5&10&14 & 4 & (9,0,2, 0,0,3, 2,0,4,0,3,0,1,  0, 0)& \frac{1}{315} \\
40a&&4&4 &14&&\{(4,10),(3,8),*\} \succ \{(7,18),* \} & \frac{1}{2520}\\
40b&&4&5 &14&&\{(2,5),(6,16),*\} \succ \{(8,21),* \} & \frac{1}{1260}\\
41& ( 0, 0, 0, 1, 1, 1, 0, 1, 1 ) & 5 & 11& 13& 2& (5,0, 1, 0,0,
0,2,0, 1, 0, 2,0, 0, 0, 0) & \frac{1}{252} \\
42&  (0, 0, 0, 1, 1, 1, 0, 1, 1) &3&6&14 & 3 &   (6,1,0, 1,0,1, 2,0,3,0,2,0,1, 0, 0) & \frac{1}{770} \\
43 &  (0, 0, 0, 1, 1, 1, 0, 1, 1) &4&8&14 & 3 &   (7,0,1, 1,0,1,1,1,2,0,2,0,1, 0, 0) & \frac{71}{27720} \\
43a&&2&4&14&& \{(4,11),(1,3),*\} \succ\{(5,14),* \} & \frac{1}{2520}\\
43b&&3&4&14&&\{(2,5),(3,8),*\} \succ \{(5,13),* \} & \frac{23}{36036}\\
43c&&3&5&14&& \{(7,16),(7,19),* \} & \frac{31}{31920}\\
44 & ( 0, 0, 0, 1, 1, 1, 0, 1, 1) &5&9&14 & 3 & (7,0,1, 1,0,1, 2,0,3,0,1,1,0, 0, 0) & \frac{43}{13860} \\
44a&&4&4&14&& \{(2,5),(6,16),*\} \succ \{(8,21),* \} & \frac{1}{1386}\\
44b&&3&3&14&& \{(7,16),(5,13),* \} & \frac{3}{16016}\\
44c&&4&6&14&& \{(7,16),(5,18),* \} & \frac{1}{720}\\
44d&&4&4&14&& \{(5,13),(5,18),* \} & \frac{1}{2184}\\
45&  (0, 0, 0, 1, 1, 1, 1, 0, 1) &4&7&14& 2 & (3,0,2, 0,0,0, 1,0,3,0,1,0,1, 0, 0) & \frac{1}{504} \\
46& ( 0, 0, 0, 1, 1, 1, 1, 1, 0) &4&7&14& 3 &( 6,0,2, 0,0,2, 1,0,3,1,1,0,1, 0, 0)& \frac{1}{504} \\
46a&&3&3&16&&\{(2,5),(3,8),*\} \succ \{(5,13),* \} & \frac{1}{16380}\\
46b&&3&6&14&& \{(3,10),(2,7),*\} \succ \{(5,17),* \} & \frac{7}{6120}\\
47&  0, 0, 0, 1, 1, 1, 1, 1, 1) &2&3&16&  2 &  (3,1,0, 1,0,0, 1,0,2,1,0,0,1, 0,0 )& \frac{1}{9240} \\
48 &  0, 0, 0, 1, 1, 1, 1, 1, 1) &3&5&14&  2 & (4,0,1, 1,0,0, 0,1,1,1,0,0,1, 0,0 )& \frac{19}{13860} \\
48a&&2&3 &14&&\{(4,9),(3,7),*\} \succ \{(7,16),* \} & \frac{1}{2640}\\
49& (  0, 0, 0, 1, 1, 1, 1, 1, 1 ) &4&7&14& 3 & (5,1,1, 0,0,2, 1,0,4, 0,2,0,1,0,0) &\frac{47}{27720} \\
49a&&4&6&14&&\{(5,11),(4,9),*\} \succ \{(9,20),* \} & \frac{1}{840}\\
50& (  0, 0, 0, 1, 1, 1, 1, 1, 1 ) &5&9&14& 3 & (6,0,2, 0,0,2, 0,1,3, 0,2,0,1,0,0) &\frac{41}{13860} \\
50a&&3&5&14&&\{(4,11),(1,3),*\} \succ \{(5,14),* \} & \frac{1}{1260}\\
51& (  0, 0, 0, 1, 1, 1, 1, 1, 1 ) &6&10&14& 3 & (6,0,2, 0,0,2, 1,0,4, 0,1,1,0,0,0) &\frac{97}{27720} \\
51a&&5&4&14&&\{(4,10),(3,8),*\} \succ \{(7,18),* \} & \frac{1}{1386}\\
51b&&5&5&14&& \{(5,13),(5,18),* \} & \frac{1}{1170}\\

52&(0, 0, 1, 0, 0, 1, 0, 1, 0 ) &3&7&14& 2 & (4,0,0, 1,0,2, 2,0,2, 0,0,0,0,0, 1) & \frac{1}{420} \\
52a&&2&3&18&&\{(2,5),(3,8),*\} \succ \{(5,13),* \} & \frac{1}{2184}\\
53&(0, 0, 1, 0, 0, 1, 1, 1, 0 ) &4&8&14& 2 & (3,0,1, 0,0,3, 1,0,3, 0,0,0,0, 0, 1) & \frac{1}{360} \\
53a&&3&4&15&&\{(2,5),(3,8),*\} \succ \{(5,13),* \} & \frac{1}{1170}\\

54& (  0, 0, 1, 0, 1, 0, 0, 1, 0) &2&4&14& 2 &( 2,0,0, 2,0,3, 1,0,1, 0,1,0,0, 0, 0)&  \frac{1}{840} \\
55 & ( 0, 0, 1, 0, 1, 0, 0, 1, 0) &2&2&14& 3 & (4,0,0, 3,0,4, 1,0,3, 0,0,1,0, 0, 0) & \frac{1}{3080}\\
56 & ( 0, 0, 1, 0, 1, 0, 1, 1, 0 ) &3&5&14& 2 & (1,0,1, 1,0,4, 0,0,2, 0,1,0,0, 0, 0)& \frac{ 1}{630} \\
56a&&2&3&14&&\{(4,9),(3,7),*\} \succ \{(7,16),* \} & \frac{1}{1680}\\
\end{array}
$$

\newpage
$$
\begin{array}{lccccccc}
No. & (P_3,\cdots,P_{11}) &P_{18}&P_{24}&\mu_1& \chi & (n_{1,2},n_{4,9},\cdots,n_{1,5})\text{ or } B_{min} & K^3\\
\hline

57& ( 0, 0, 1, 0, 1, 0, 1, 1, 0 ) &3&3&14& 3 & (3,0,1, 2,0,5, 0,0,4, 0,0,1,0, 0, 0) & \frac{1}{1386} \\
58& (0, 0, 1, 0, 1, 1, 0, 1, 0) &3&6&14& 3 & (4,0,1, 1,0,4, 2,0, 2, 0,1,0,1, 0, 0) & \frac{1}{630} \\
58a&&2&4&14&&\{(4,9),(3,7),*\} \succ \{(7,16),* \} & \frac{1}{1680}\\
59& (0, 0, 1, 0, 1, 1, 0, 1, 1) &2&4&14& 2 & (2,0,0, 2,0,2, 1,1,0, 0,0,0,1, 0, 0) &  \frac{3}{3080} \\
59a&&2&3&14&&\{(3,8),(4,11),*\} \succ \{(7,19),* \} & \frac{1}{2660}\\
60& (0, 0, 1, 0, 1, 1, 1, 1, 0) &4&7&14& 3 & (3,0, 2, 0,0,5, 1,0,3, 0,1,0,1, 0, 0 ) & \frac{1}{504} \\
60a&&3&3&15&&\{(2,5),(3,8),*\} \succ \{(5,13),* \} & \frac{1}{16380}\\
61& (0, 0, 1, 0, 1, 1, 1, 1, 1 ) &2&3&15& 2 & (0,1,0, 1,0, 3,1,0,2,0,0,0,1,0,0 )& \frac{1}{9240} \\
62& (0, 0, 1, 0, 1, 1, 1, 1, 1 ) &3&5&14& 2 & (1,0,1, 1,0,3, 0,1,1,0,0,0,1,0,0 )& \frac{19}{13860} \\
62a&&2&3&14&&\{(4,9),(3,7),*\} \succ \{(7,16),* \} & \frac{1}{2640}\\

63&(0, 0, 1, 1, 1, 1, 1, 1, 1 ) &3&4&14& 3 & (5,0,1, 2,0,1, 1,1,3, 0,1,0,0,0, 1) & \frac{1}{5544} \\
\end{array}
$$}





%

\begin{proof} Assume $B=B_{*}$. We hope to deduce a contradiction.

(1). If $P_{14} \ge 2$, then Proposition \ref{volume} implies $K^3
\ge \frac{11}{37800}>\frac{1}{3437}$. Thus $B\neq B_{4a}$, $B_{9}$,
$B_{16c}$, $B_{24}$, $B_{27a}$, $B_{29a}$, $B_{44b}$, $B_{63}$.

(2). If $P_{15} \ge 2$, then Proposition \ref{volume} implies $K^3
\ge \frac{11}{46080}
> \frac{1}{4190}$.  Hence $B\neq B_{60a}, B_{61}$.

(3). If $P_{16} \ge 2$, then Proposition \ref{volume} implies $K^3
\ge \frac{11}{55488} > \frac{1}{5045}$. Hence $B\neq  B_{46a},
B_{47}$.

(4). If $P_{18} \ge 2$, then Proposition \ref{volume} implies $K^3
\ge \frac{11}{77976}>\frac{1}{7089}$. Thus $B\neq  B_{20a},
B_{22}$.

(5).  Besides, we see $P_6(B_{33a})=1$, $P_{16}(B_{33a})=2$ but
$P_{22}(B_{33a})=1$, a contradiction. So $B\neq B_{33a}$.

(6). For cases  16a, 18a, 21a, 52a and case 55, one has
$P_{17}(B_{*})=0$. But since $P_8(B_{21a})=P_9(B_{21a})=1$, $B\neq
(B_{21a})$. Also for case 52a and case 55, since
$P_5(B_{*})=P_{12}(B_{*})=1$, we see $B\neq B_{52a}, B_{55}$. For
case 18a, since $P_6(B_{18a})=P_{11}(B_{18a})=1$, we see $B\neq
B_{18a}$. Finally since  $P_{19}(B_{16a})=-1$, we see $B\neq
B_{16a}$.
\end{proof}

\begin{thm}\label{vol2} Assume $\chi(\OO_X)>1$. Then
$K_X^3 \geq \frac{1}{2660}$. Furthermore, $K_X^3= \frac{1}{2660}$
if, and only if, $P_2=0$ and either $\chi=3$, $B=B_{11a}$ or
$\chi=2$, $B_{59a}$.
\end{thm}
\begin{proof} If $P_{m_0} \geq 2$ for some positive integer $m_0
\le 12$, then Proposition \ref{volume} implies $K_X^3 \geq
\frac{11}{24336}>\frac{1}{2213}>\frac{1}{2660}$.

Assume $P_m \leq 1$ for $m \leq 12$. Then we have seen $B\geq B_{*}$
where $B_{*}$ is one in Table C excluding those cases listed in
Proposition \ref{impos}.

We can see $K^3(B_{11a})=K^3(B_{59a})=\frac{1}{2660}$.

We pick out those cases with $K^3(B_{*})<\frac{1}{2660}$. They are
cases 4a, 16a, 16c, 18a, 20a, 21a, 27a, 29a, 33a, 44b, 46a and case
60a. In all these cases, Corollary \ref{impos} says $B\neq B_{*}$.
Thus $B\succ B_{*}$. In order to prove the theorem, we need to
study the one step unpacking of $B_{*}$ case by case.


%


First we consider case 4a and case 29a. It's obtained by 2-steps of
packing from $B_4$:
$$B_4=\{(2,6),(4,11),*\} \succ B_{4.5}:=\{(1,3),(5,14),*\} \succ \{(6,17),*\}=B_{4a}.$$
By \cite[Lemma 3.6]{Explicit_I}, we get $K_X^3=K^3(B)\geq
K^3(B_{4.5})=\frac{1}{630}> \frac{1}{2660}$. Similarly, we also get
$K_X^3>\frac{1}{2660}$ for case 29a.

Next we consider cases 18a, 20a, 21a, 46a, 52a, 60a. The common
property is that they are obtained by a 1-step packing from
$B^{(12)}$. So the only possibility is $B^{(12)}=B$. Thus
$K_X^3=K^3(B_{18})$ or $K^3(B_{20})$ or $K^3(B_{21})$ or
$K^3(B_{46})$ or $K^3(B_{52})$ or $K^3(B_{60})$. In a word,
$K_X^3>\frac{1}{2660}$.

The remaining cases are: 16a, 16c, 27a and 44b. For case 44b, there are two
intermediate baskets dominating $B_{44c}$ or $B_{44d}$,
respectively. Thus, in particular, $K_X^3> \frac{1}{2184}$. For case
27a, it's obtained from $B_{27}$ by 3-steps of packing, namely
$$\begin{array}{lll}
&&B_{27}=\{3 \times (2,5), (5,8),* \} \succ \{ 2 \times (2,5),
(5,13),* \}\\
& \succ &B_{27.5}:=\{ (2,5),(7,18),* \} \succ \{
(9,23),*\}=B_{27a}.\end{array}$$ Thus we see: $B\succcurlyeq
B_{27.5}$ and $K_X^3\geq K^3(B_{27.5})=
\frac{1}{1386}>\frac{1}{2660}$. Finally we consider cases 16a and
16c. We know $$B_{16}=\{(4,9),(3,7),(2,5),(2,5),(3,8),(3,8),*\}.$$
The 1-step packing of $B_{16}$ yields
$$B_{16.5}:=\{(4,9),(3,7),(2,5),(5,13),(3,8),*\}$$ and the 1-step
prime packing of $B_{16.5}$ is either $B_{16a}$ or $B_{16c}$. Thus,
if $B\succcurlyeq B_{16.5}$, then $K_X^3\geq K^3(B_{16.5})=
\frac{85}{72072}>\frac{1}{848}$. The other intermediate basket
dominating $B_{16a}$ and $B_{16c}$ is
$$B_{16.6}:=\{(7,16),(2,5),(2,5),(3,8),(3,8),*\}$$ with
$K^3(B_{16.6})=\frac{13}{6160}>\frac{1}{474}$. No other ways to
obtain either $B_{16a}$ or $B_{16c}$ beginning from $B^{(16)}$. The
theorem is proved.
\end{proof}

\begin{cor} \label{p24} Assume $\chi(\OO_X)>1$. Then
$P_{24}\geq 2$.
\end{cor}
\begin{proof} By \cite[Theorem 4.15]{Explicit_I}, we know either
$P_{10}\geq 2$ or $P_{24}\geq 2$. When $q(X)>0$, the statement
follows from \cite{JC-H}. So we may assume $q(X)=0$.

If $P_{10}\geq 2$, we take $m_0=10$ and study $\varphi_{10}$. Keep
the same notation as in \ref{setup}. By Lemma \ref{chi(O)>1}, $f$
must be of type $III$, $II$, $I_p$. Proposition \ref{nonvanishing}
(i), Theorem \ref{tIII} (1) and Theorem \ref{tII} (1) imply
$P_{24}\geq 2$.
\end{proof}


\begin{thm}\label{18} Assume $\chi(\OO_X)>1$.
Then $P_{m_0}\geq 2$ for certain positive integer $m_0\leq 18$. In
particular, $\mu_1\leq 18$.
\end{thm}

\begin{proof}
Assume $P_{m} \leq 1$ for all $m\leq 12$. Then Table C tells that
$$B^{(12)}\succcurlyeq B\succcurlyeq B_{\text{min}}$$ where $
B_{\text{min}}$ is of certain type in Table C. Since, in Table C,
we have seen $\mu_1(B_{\text{min}})\leq 18$, thus \cite[Lemma
3.6]{Explicit_I} implies $\mu_1(X)\leq \mu_1(B_{\text{min}})\leq
18$.
\end{proof}

\begin{thm}\label{rho} Assume $\chi(\OO_X)>1$. Then $\rho_0(X)\leq 27$.
\end{thm}

\begin{proof} The statement follows from \cite{JC-H} when $q(X)>0$.
Assume $q(X)=0$ from now on.

If $P_{m_0} \ge 2$ for some $m_0 \le 12$, then the induced
fibration $f$ from $\varphi_{m_0}$ is of type $III$, $II$ or $I_p$
by Lemma \ref{chi(O)>1}. Thus Proposition \ref{nonvanishing} (i),
Theorem \ref{tIII} (1) and Theorem \ref{tII} (1) imply that $P_m>0$
for all $m\geq 27$.

If $P_{m}\leq 1$ for all $m\leq 12$, we have a complete
classification (cf. Table C). For each $B_{\text{min}}$ in Table C,
we observed that $P_m>0$ for all $47 \geq m\geq 24$. This is enough
to assert $P_m>0$ for all $m\geq 24$. We are done.
\end{proof}

\section{\bf Pluricanonical birationality}

In this section, we mainly study the birationality of $\varphi_m$.
Then we can conclude our main theorems. Let $X$ be a projective
minimal 3-fold of general type. First, we recall several known
theorems.

\begin{thm}\label{q>0}  (\cite{JC-H}) Assume $q(X):=h^1(\OO_X)>0$. Then $\varphi_m$
is birational for all $m\geq 7$.
\end{thm}

\begin{thm}\label{5k+6} (\cite[Theorem 0.1]{JPAA}) Assume $P_{m_0}\geq 2$ for some positive integer $m_0$.
Then $\varphi_{m}$ is birational onto its image for all $m\geq
5m_0+6$.
\end{thm}

\begin{thm}\label{CZ} (\cite{Chen-Zuo}) Assume $\chi(\OO_X)\leq 0$.
Then $\varphi_m$ is birational for all $m\geq 14$.
\end{thm}

We need the following lemma to prove our main theorems.

\begin{lem}\label{I-a(2)(3)} Assume $P_{m_0}(X)\geq 2$ for some positive integer
$m_0$. Keep the same notation as in \ref{setup} and assume $f$ is
of type $I_p$ or $I_n$. Suppose $|G|$ is a base point free linear
system on $S$. If there exists an integer $m_1>0$ with
$m_1\pi^*(K_X)|_S\geq G$, then Assumptions \ref{asum} (2),
 is  satisfied for all integers $$m\geq
\text{max}\{\rho_0+m_0+m_1, m_0+m_1+2\}.$$
\end{lem}
\begin{proof}
Since $$\begin{array}{rrl}
&& K_S+\roundup{(m-1)\pi^*(K_X)-S-\frac{1}{p}E_{m_0}'}_{|S}\\
&\geq & K_S+(m-1)\pi^*(K_X)_{|S}-(S+E_{m_0}')_{|S}\\
&\geq & (m-m_0)\pi^*(K_X)_{|S} \geq (m-m_0-m_1)\pi^*(K_X)_{|S}+G
\end{array}$$
and
$$\begin{array}{lll}
&&K_S+(m-m_0-1)\pi^*(K_X)_{|S}\\
&\geq& K_S+(m-m_0-m_1-1)\pi^*(K_X)_{|S}+G,
\end{array}$$
Lemma \ref{S2} implies that
$|K_S+\roundup{(m-1)\pi^*(K_X)-S-\frac{1}{p}E_{m_0}'}_{|S}|$ can
distinguish different generic irreducible elements of $|G|$ when
$m\geq \rho_0+m_0+m_1$ and $m\geq m_0+m_1+2$.
\end{proof}

\begin{thm}\label{chi-1}\footnote{L. Zhu \cite{Zhu2} showed $\varphi_m$ is birational for $m\geq 46$.} Let $X$ be a projective minimal 3-fold of
general type with $\chi(\OO_X)=1$. Then $\varphi_m$ is birational
for all $m\geq 40$.
\end{thm}
\begin{proof} If $P_{m_0}\geq 2$ for some $m_0\leq 6$, then, by
Theorem \ref{5k+6}, $\varphi_m$ is birational for $m\geq 36$.

Assume $P_m\leq 1$ for all $m\leq 6$. Then, by Corollary \ref{min},
$\mathscr{B}(X)$ either dominates a minimal basket in $$\{B_{2,1},
B_{2,2}, B_{3,1}\sim B_{3,5}, B_{5,1}\sim B_{5,3}, B_{6,1}\sim
B_{6,6}, B_{8,1}, B_{8,2}\}$$ or dominates the basket $B_{210}$. We
have known $P_m(X)\geq P_m(B_{*,*})$. By analyzing all the above
baskets, we see a common property that there is a pair of positive
integers $(n_0,n_1)$ satisfying $P_{n_0}\geq 2$, $P_{n_1}\geq 3$ and
one of the following conditions:

(1) $n_0\leq 10$, $n_1\leq 12$ (see cases III-2, III-3, III-4,
VI-6);

(2) $n_0\leq 9$, $n_1\leq 13$ (for the rest cases).

By Theorem \ref{7} (1), we know $\rho_0\leq 7$. We set $m_0=n_1$.
Keep the same notation as in \ref{setup}. Our proof is organized
according to the type of $f$. Note that $P_{m_0}\geq 3$ and
$m_0\leq 13$. By Theorem \ref{q>0}, we only need to care about the
situation: $q(X)=0$.
\medskip

%

{\bf Case 1}. $f$ is of type $I_3$.

Take $G$ to be the movable part of $|2\sigma^*(K_{S_0})|$. Claim A
implies that Assumptions \ref{asum} (2) is satisfied whenever
$m\geq 35\geq \rho_0+2m_0+2$. Clearly, by Lemma \ref{a(1)},
Assumptions \ref{asum} (1) is also satisfied. As seen in the later
part of \ref{I3}, we can take a rational number $\beta\mapsto
\frac{p}{2m_0+2p}\geq \frac{1}{m_0+2}$.  Now inequality (2.2) gives
$\xi\geq \frac{4}{45}$. Take $m=35$. Then
$\alpha=(35-1-\frac{m_0}{2}-\frac{1}{\beta})\xi\geq
\frac{10}{9}>1$. Theorem \ref{technical} gives $\xi\geq
\frac{4}{35}$. Take $m=32$. Then $\alpha>1$. We will see $\xi\geq
\frac{1}{8}$ similarly. Now, for $m\geq 39$, $\alpha\geq
(39-1-\frac{13}{2}-15)\xi\geq \frac{33}{16}>2$. Theorem
\ref{technical} says that $\varphi_m$ is birational for all $m\geq
39$.
\medskip

{\bf Case 2}. $f$ is of type $II$ or $III$.

We take $\tilde{m}_0=n_0$ and $m_1=n_1$. We still use the mechanics
of \ref{setup} to study $\varphi_{\tilde{m}_0}$ in stead of
$\varphi_{m_0}$. But most notations will be put on the symbol
$\widetilde{\ \ }$. Noting that $\tilde{m}_0\leq 10$ and
$P_{\tilde{m}_0}\geq 2$.

If $\tilde{f}$ is of type $II$ or $III$, Theorem \ref{tIII} (3) and
Theorem \ref{tII} (3) imply that $\varphi_m$ is birational for
$m\geq 38$.

If $\tilde{f}$ is of type $I_n$ or $I_p$. We take $\tilde{G}$ to be
the movable part of $|M_{m_1}|_{\tilde{S}}|$, where $\tilde{S}$ is a
generic irreducible element of $|M_{\tilde{m}_0}|$. Clearly
$h^0(\tilde{S}, M_{m_1}|_{\tilde{S}})\geq 2$ since
$\dim\varphi_{m_1}(X)\geq 2$. Thus we are in the situation with
$m_1\pi^*(K_X)|_{\tilde{S}}\geq \tilde{G}$. We may always take a
sufficiently good $\tilde{\pi}$ instead of $\pi$. Now Lemma
\ref{a(1)} and Lemma \ref{I-a(2)(3)} imply that Assumptions
\ref{asum} (1), (2) are simultaneously satisfied for $m\geq 30\geq
\rho_0+\tilde{m}_0+m_1$. {}Finally, we study the value of $\alpha$.
Clearly, one may take $\tilde{\beta}=\frac{1}{m_1}$. Thus inequality
(2.2) says $\xi\geq \frac{2}{1+\tilde{m_0}+m_1}$. For situations (1)
and (2), we have $\xi\geq \frac{2}{23}$. Take $m=35$. Then
$\alpha\geq \frac{24}{23}>1$. Theorem \ref{technical} gives $\xi\geq
\frac{4}{35}$. Take $m=32$. Then similarly we get $\xi\geq
\frac{1}{8}$. Take $m\geq 40$. Then $\alpha\geq \frac{17}{8}>2$.
Theorem \ref{technical} implies that $\varphi_m$ is birational for
all $m\geq 40$. We are done.
\end{proof}

\begin{thm}\label{birat}
Let $X$ be a projective minimal 3-fold of general type with
$\chi(\OO_X)>1$. Then $\varphi_m$ is birational for all $m \ge 73$.
\end{thm}

\begin{proof} By Theorem \ref{q>0}, we only need to consider the
situation $q(X)=0$. According to Lemma \ref{chi(O)>1}, the induced
fibration $f$ from $\varphi_{m_0}$ is of type $III$, $II$ or $I_p$.

If $P_{m_0} \ge 2$ for some $m_0 \le 16$, then, by Theorems
\ref{tIII}, \ref{tII} and \ref{tIp}, $\varphi_m$ is birational for
all $m \ge 69$ . Assume $P_m\leq 1$ for all $m\leq 16$. Then we have
a complete classification for $B_{\text{min}}$ as in Table C. More
precisely, we see $B\succcurlyeq B_{2a}$, $B\succcurlyeq B_{3a}$ and
$B\succ B_{20a}$, $B\succ B_{52a}$, noting that case No.22 doesn't
happen by Proposition \ref{impos}. As we have observed in the proof
of Theorem \ref{vol2}, for cases No. 20a and No. 52a, we actually
have $B=B_{20}$ and $B=B_{52}$. Thus we see $P_{14}(X)\geq 2$ in
both cases, a contradiction. We are left to study cases: No. 2a and
No. 3a, which correspond to two formal baskets: $(B_{2a},2,0)$ and
$(B_{3a},3,0)$, where
$$B_{2a}=\{4\times (1,2), (4,9), (2,5), (5,13), 3\times (1,3),
2\times (1,4)\}$$
$$B_{3a}=\{6\times (1,2), (5,11), 3\times (2,5), (5,13), 4\times
(1,3), (2,7), 2\times (1,4)\}.$$ The computation gives the following
datum:
\medskip

\begin{center}
\begin{tabular}{|c|c|c|c|c|}
\hline
&$\rho_0$&$\mu_1$&$\mu_2$&$\mu_3$\\
\hline
$B_{2a}$&20&18&24&30\\
\hline $B_{3a}$&20&18&20&30\\
 \hline
\end{tabular}\end{center}

When $B\succcurlyeq B_{3a}$, we have $P_{20}(X)=P_{20}(B)\geq
P_{20}(B_{3a})\geq 3$. Theorems \ref{tIII} and \ref{tI3}  imply that
$\varphi_m$ is birational for all $m\geq 66$ unless $f$ is type
$II$. Indeed, if $f$ is of type $II$ and $m_0=20$, at least we have
$\xi \ge \frac{2}{31}$ following the argument in \ref{II}. Take
$m=57$, we have $\alpha
>1$ and hence $\xi \ge \frac{4}{57}$. Now take $m\geq 70$, we have
$\alpha=(70-41)\frac{4}{57}>2$. Thus $\varphi_m$ is birational for
all $m \ge 70$.

 Now the theorem follows from the following claim.
\medskip

{\bf Claim B.} When $B\succcurlyeq B_{2a}$, $\varphi_m$ is
birational for all $m\geq 73$.
\medskip

The proof is similar to that of Theorem \ref{chi-1}, Case 1 and Case
2. We have known $\rho_0\leq 20$. We can find two numbers $n_0\leq
18$ and $n_1\leq 24$ with $P_{n_0}(X)\geq 2$ and $P_{n_1}(X)\geq 3$.
{}First, we set $m_0=n_1$. Keep the same notation as in \ref{setup}.
Our proof is organized according to the type of $f$. Note that
$P_{m_0}\geq 3$ and $m_0\leq 24$.
\medskip

{\bf Case i}. $f$ is of type $I_3$.

By Lemma \ref{chi(O)>1}, $f$ must be of type $I_p$. Take
$G=2\sigma^*(K_{S_0})$. Claim A implies that Assumptions \ref{asum}
(2) is satisfied whenever $m\geq 70\geq \rho_0+2m_0+2$. Clearly, by
Lemma \ref{a(1)}, Assumptions \ref{asum} (1) is also satisfied. As
seen in the latter part \ref{I3}, we can take a rational number
$\beta\mapsto \frac{p}{2m_0+2p}\geq \frac{1}{m_0+2}$.  Note that
$|G|$ is base point free, we have $\deg(K_C)\geq 6$. Now inequality
(2.2) gives $\xi\geq \frac{2}{13}$. For $m\geq 70$, $\alpha\geq
(70-1-12-26)\xi >2$. Theorem \ref{technical} says that $\varphi_m$
is birational for all $m\geq 70$.
\medskip

{\bf Case ii}. $f$ is of type $II$ or $III$.

We take $\tilde{m}_0=n_0$ and $m_1=n_1$. We still use the mechanics
of \ref{setup} to study $\varphi_{\tilde{m}_0}$ in stead of
$\varphi_{m_0}$. Noting that $\tilde{m}_0\leq 18$, when $\tilde{f}$
is of type $III$ or $II$, Theorems \ref{tIII} and \ref{tII} imply
that $\varphi_m$ is birational for all $m\geq 66$. We are left to
study the situation with $\tilde{f}$ being of type $I$. We take
$\tilde{G}$ to be the movable part of $|M_{m_1}|_{\tilde{S}}|$.
Clearly $h^0(\tilde{S}, M_{m_1}|_{\tilde{S}})\geq 2$ since
$\dim\varphi_{m_1}(X)\geq 2$. Thus we are in the situation with
$m_1\pi^*(K_X)|_{\tilde{S}}\geq \tilde{G}$. Now Lemma \ref{a(1)} and
Lemma \ref{I-a(2)(3)} imply that Assumptions \ref{asum} (1), (2) are
simultaneously satisfied for $m\geq 62\geq \rho_0+\tilde{m}_0+m_1$.
Clearly, one may take $\tilde{\beta}=\frac{1}{m_1}$. Thus inequality
(2.2) says $\xi\geq \frac{2}{1+\tilde{m_0}+m_1}\geq \frac{2}{43}$.

Take $m=65$. Then $\alpha\geq \frac{44}{43}>1$. Theorem
\ref{technical} gives $\xi\geq \frac{4}{65}$. Take $m=60$. Then
similarly we get $\xi\geq \frac{1}{15}$. Take $m=59$. Then we shall
get $\xi\geq \frac{4}{59}$. Take $m=58$ and we obtain $\xi\geq
\frac{2}{29}$. Eventually, for $m\geq 73$, we see $\alpha>2$ and
Theorem \ref{technical} implies that $\varphi_m$ is birational for
all $m\geq 73$. We are done.
\end{proof}

We have proved all the main results. Indeed, Theorem
\ref{1} follows from Theorem \ref{q>0}, Theorem \ref{CZ}, Theorem
\ref{420} and Theorem \ref{chi-1}. Theorem \ref{m} follows from
Theorem \ref{1}, Theorem \ref{vol2}, Corollary \ref{p24}, Theorem
\ref{18}, Theorem \ref{rho}, Theorem \ref{birat}.

We would like to propose the
following:

\begin{conj} For all projective minimal 3-folds of general type, the
inequality $K^3\geq \frac{1}{2660}$ is optimal.
\end{conj}

The following problem is  very interesting.

\begin{op} Can one find a minimal 3-fold $X$ of general type with $q(X)=0$ and
$\chi(\OO_X)>1$?
\end{op}


\end{document}